\documentclass[opre,nonblindrev]{informs3}

\DoubleSpacedXI

\usepackage[utf8]{inputenc}
\usepackage{booktabs}
\usepackage{multirow}
\usepackage{enumerate}
\usepackage{natbib}
\usepackage{graphicx}
\usepackage{tikz}
\usepackage{subcaption}
\usepackage{amssymb}
\usepackage{array}
\usepackage{float}
\usepackage{xcolor}
\usepackage{comment}
\usepackage[nocomma,short]{optidef} 
\usepackage{bm}
\usepackage{mathtools}
\usepackage{pgfplots}

\definecolor{gray1}{gray}{1}
\definecolor{gray2}{gray}{0.95}
\definecolor{gray3}{gray}{0.8}
\definecolor{gray4}{gray}{0.65}

\usepackage{hyperref}
\usetikzlibrary{positioning}
\hypersetup{colorlinks=true,linkcolor=black,citecolor=blue,filecolor=black,urlcolor=black,pdfpagelayout=SinglePage}
\newcolumntype{C}[1]{>{\centering\arraybackslash}p{#1}}
\usetikzlibrary{patterns}
\usetikzlibrary{arrows,automata}
  \usetikzlibrary{shapes,arrows}
  \usetikzlibrary{shapes.geometric}
  
  \usetikzlibrary{arrows,decorations.pathmorphing,backgrounds,positioning,fit,shapes.geometric,calc,shapes}
\usetikzlibrary{shapes,shapes.geometric,arrows,fit,calc,positioning,automata}
\tikzset{triangle/.style={regular polygon, regular polygon sides=3}}
\tikzstyle{arc} = [draw, very thick, color=black, -latex']
\tikzstyle{v1} = [draw, very thick, color=blue, -latex']
\tikzstyle{v2} = [draw, very thick, color=red, -latex']

\tikzstyle{t1} = [draw, very thick, color=red, -latex']
\tikzstyle{t2} = [draw, very thick, color=red, -latex']

\tikzstyle{t1v1} = [draw, dashed, color=blue]
\tikzstyle{t2v1} = [draw, dashed, color=blue]
  
\usetikzlibrary{shapes,arrows}
\usetikzlibrary{shapes.geometric}
\tikzstyle{decision} = [diamond, draw,fill=white,
    text width=4em,
    text badly centered, inner sep=0pt]
\tikzstyle{block} = [rectangle, draw, fill=white,
    text width=8em,
    text centered,
    minimum height=4em,
    minimum width = 6em
    ]
\tikzstyle{line} = [draw, very thick, color=black, -latex']
\tikzstyle{cloud} = [draw, ellipse,fill=red!20, node distance=2.5cm,minimum height=2em
    ]
\usepackage{eurosym}
\usepackage{mathtools}
\bibpunct[, ]{(}{)}{,}{a}{}{,}%
\def\bibfont{\small}%

\usepackage{times}
\usepackage{booktabs}
\usepackage{enumerate}
\usepackage{mathrsfs}
\usepackage{hyperref}
\PassOptionsToPackage{hyphens}{url}
\usepackage{xurl}
\usepackage{mathtools}
\usepackage{algorithmicx}
\usepackage{algcompatible}
\usepackage{pdflscape}
\usepackage{multicol}
\usepackage[english]{babel}
\usepackage[utf8]{inputenc}
\hypersetup{colorlinks=true,linkcolor=black,citecolor=blue,filecolor=black,urlcolor=blue,pdfpagelayout=SinglePage, linkbordercolor=blue,pdfborderstyle={/S/U/W 1}}
\usepackage{esvect}
\usepackage{enumitem}
\usepackage{pgf}
 \pgfplotsset{compat=1.17} 
\usetikzlibrary{arrows,decorations.pathmorphing,backgrounds,positioning,fit,shapes.geometric,calc,shapes}

\usepackage[norelsize, linesnumbered, ruled, lined, boxed, commentsnumbered]{algorithm2e}
\usepackage{amsmath, amsfonts}
\usepackage{array,delarray}   
\usepackage{makecell}

\usepackage{mathbbol}

\newcommand{\st}{\text{s.t.}}

\newcommand{\mcl}[1]{\mathcal{#1}}
\newcommand{\mbb}[1]{\mathbb{#1}}
\newcommand{\ov}[1]{\overline{#1}}
\newcommand{\un}[1]{\underline{#1}}

\newcommand{\til}[1]{\tilde{#1}}

\newcommand{\bbb}{\{0,1\}}

\newcommand{\rrr}{\mbb{R}}

\newcommand{\snode}{s}
\newcommand{\dnode}{d}

\newcommand{\travt}[2]{\bm{t}_{{#1}{#2}}}

\newcommand{\Kset}{\mathcal{K}}

\newcommand{\uncSet}{\Xi}

\newcommand{\alDV}{\boldsymbol{\alpha}}

\newcommand{\betaDV}{\boldsymbol{\beta}}
\newcommand{\betakDV}[1]{\boldsymbol{\beta}^{#1}}
\newcommand{\gammaDV}[1]{\boldsymbol{\gamma}^{#1}}

\newcommand{\ySet}{\mathcal{Y}}
\newcommand{\discV}{\bm{w}}
\newcommand{\discVdom}{\mathcal{W}}
\newcommand{\yvk}[1]{\bm{y}^{#1}}
\newcommand{\tyvk}[1]{\til{\bm{y}}^{#1}}
\newcommand{\xvk}[1]{\bm{z}^{#1}}

\newcommand{\uncV}{\boldsymbol{\xi}}

\newcommand{\unckV}[1]{{\boldsymbol{\xi}^{#1}}}

\newcommand{\edgeSet}{\mcl{E}}

\newcommand{\nodeSet}{\mcl{N}}

\newcommand{\cutset}[1]{\mcl{C}({#1})}

\newcommand{\sol}{\left(\discV, \alDV, \betaDV, \{\yvk{k},\betakDV{k},\gammaDV{k}\}_{k \in \Kset}\right)}

\newtheorem{observation}{Observation}

\MANUSCRIPTNO{}

\usepackage{endnotes}
\let\footnote=\endnote

\usepackage{natbib}
 \bibpunct[, ]{(}{)}{,}{a}{}{,}%
 \def\bibfont{\small}%
\TheoremsNumberedThrough     
\ECRepeatTheorems

\EquationsNumberedThrough    

\begin{document}

\RUNAUTHOR{Paradiso, Georghiou, Dabia, Tönissen}

\RUNTITLE{Exact and Approximate Schemes for Problems with Decision Dependent Information Discovery}

\TITLE{\vspace{-2cm}~\\Exact and Approximate Schemes\\ for Robust Optimization Problems\\ with Decision Dependent Information Discovery}

\ARTICLEAUTHORS{%
\AUTHOR{Rosario Paradiso}
\AFF{Department of Operations Analytics, Vrije Universiteit Amsterdam, the Netherlands, \EMAIL{r.paradiso@vu.nl}} 
\AUTHOR{Angelos Georghiou}
\AFF{Department of Business and Public Administration, University of Cyprus, Cyprus, \EMAIL{georghiou.angelos@ucy.ac.cy}}
\AUTHOR{Said Dabia}
\AFF{Department of Operations Analytics, Vrije Universiteit Amsterdam, the Netherlands, \EMAIL{s.dabia@vu.nl}}
\AUTHOR{Denise Tönissen}
\AFF{Department of Operations Analytics, Vrije Universiteit Amsterdam, the Netherlands, \EMAIL{d.d.tonissen@vu.nl}}
} 
\ABSTRACT{
\noindent Uncertain optimization problems with decision dependent information discovery allow the decision maker to control the timing of information discovery, in contrast to the classic multistage setting where uncertain parameters are revealed sequentially based on a  prescribed filtration. This problem class is useful in a wide range of applications, however, its assimilation is partly limited by the lack of efficient solution schemes. In this paper we study two-stage robust optimization problems with decision dependent information discovery where uncertainty appears in the objective function. The contributions of the paper are twofold: $(i)$ we develop an exact solution scheme based on a nested decomposition algorithm, and $(ii)$ we improve upon the existing $K$-adaptability approximate  by strengthening its formulation using techniques from the integer programming literature. Throughout the paper we use the orienteering problem as our working example, a challenging problem from the logistics literature which naturally fits within this framework. The complex structure of the routing recourse problem forms a challenging test bed for the proposed solution schemes, in which we show that exact solution method outperforms at times the $K$-adaptability approximation, however, the strengthened $K$-adaptability formulation can provide good quality solutions in larger instances while significantly outperforming existing approximation schemes even in the decision independent information discovery setting. We leverage the effectiveness of the proposed solution schemes and the orienteering problem in a case study from Alrijne hospital in the Netherlands, where we try to improve the collection process of empty medicine delivery crates by co-optimizing sensor placement and routing decisions.

}%



\maketitle

%


\section{Introduction}
Robust optimization has proven to be a successful methodology to tackle optimization problems under uncertainty. In the last decades, rich and consolidated theory has been developed for the case in which all the decisions are taken \emph{before} uncertainty is unfolded, i.e., \emph{single stage robust optimization} with \emph{here-and-now decisions} (\cite{ben2009robust,bertsimas2011theory,gabrel2014recent}). 
When some decisions can be adjusted after observing uncertainty (\emph{recourse} or \emph{wait-and-see} decisions), the problem falls into the domain of \emph{multistage robust optimization}. These problems not only can capture the dynamic nature of the process but they can also mitigate the conservatism of the static counterpart. On the other hand, they are noticeably more challenging and still resist to solution in their general form.
A rich stream of research tries to tackle problems by approximating the functional form of recourse decisions, known as the decision rule approximation. Decision rule approximations have been developed for continuous (\cite{ben2004adjustable,chen2009uncertain,georghiou2015generalized,bertsimas2019adaptive}) as well as binary (\cite{bertsimas2015design,postek2016multistage,bertsimas2016multistage,bertsimas2018binary}) recourse decisions, see \cite{delage2015robust} for a comprehensive review.
When the problem features binary recourse decisions, one of the most successful methodologies is the $K$-adaptability approximation (\cite{bertsimas2010finite,hanasusanto2015k,buchheim2017min,chassein2019faster,subramanyam2019k}). This approximation selects $K$ candidate constant policies before observing the realization of the uncertain parameters and then implements the best of these policies after the realization is known.

 In the classic multistage setting, it is assumed that uncertain parameters are sequentially revealed to the decision-maker based on a prescribed filtration. Yet, this assumption fails to hold in a number of applications where the actions of the decision-maker can influence the timing of information discovery. For example, in \emph{clinical trials} the success of the clinical study can only be revealed after the pharmaceutical has chosen to undertake the study, in \emph{offshore gas-field developments} the quantity of gas reserves can only be observed if an exploratory drill is performed, and in \emph{production planning} the production cost of a new product is uncertain and cannot be observed until the company choose to produce it (see Table \ref{tab:literature} for more examples). This class of problems is termed as problems with \emph{decision dependent information discovery} (\cite{vayanos2011decision,vayanos2020robust}). The majority of these applications are modeled through the lens of stochastic programming where controlling the non-anticipative structure of the decisions is translated as controlling the structure of the scenario tree using binary decisions (\cite{goel2006class,gupta2011solution,boland2016minimum}). In the robust optimization setting \cite{vayanos2011decision} leverages the decision rule approximation scheme to control the timing of information discovery, while  \cite{vayanos2020robust} provide a more scalable and accurate approximation by extending the $K$-adaptability approximation to the multistage decision dependent information discovery setting.

 Many applications from the field of logistics can be cast as problems with decision dependent information discovery. Yet, to the best of our knowledge, these synergies have not been explored with the exception of \cite{macias2020endogenous} which models the monitoring and relief of a humanitarian distribution network through the use of UAVs. In this paper, we use as our working example the \emph{orienteering problem} which is a challenging logistics problem. First introduced by \cite{tsiligirides1984heuristic},  the deterministic orienteering problem is combination of the knapsack  and  the traveling salesman problem described as follows: given a set of locations each with an assigned profit, the objective is to compute the path that maximizes the total profit collected while not violating the time duration constraint. In its basic form, the orienteering problem has found many applications including distribution, maintenance, tourism, mobile crowd-sourcing and military applications, see the review papers (\cite{vansteenwegen2011orienteering,gunawan2016orienteering}) for further applications and variants of the problem. Stochastic or robust variants of the orienteering problem have been proposed in (\cite{campbell2011orienteering,evers2014robust,evers2014two,verbeeck2016solving}), however, all works consider the information discovery being independent of the actions of the decision maker. Yet, many  logistic problems could benefit from this modeling approach. For example, the tourist trips design problem (\cite{vansteenwegen2007mobile,gavalas2014survey}) can be combined with the preference elicitation problem (\cite{bertsimas2013learning,vayanos2020active}) to learn the uncertain preferences of  individual tourists thus being able to propose personalized sightseeing recommendations.  In this two-stage robust decision dependent information discovery setting, the here-and-now discovery decisions are the questions used to ``observe" the unknown ``rewards" experienced by the tourist at each attraction, while the recourse actions decide the routing plan the tourist will follow to visit the attractions  which is formulated as an orienteering problem.  An example which we leverage in this study, and which is motivated by the delivery of medicine in the Alrijne hospital in the Netherlands (\cite{alrijneMasterThesis,alrijneBachelorThesis}), tries to improve the collection  process of empty medicine delivery crates.  In this two-stage setting, the here-and-now discovery decisions are placements of sensors in different departments of the hospital to detect the uncertain number of empty crates while the recourse actions decide the collection route which we formulated as an orienteering problem. 

\begin{table}[htbp!]
\begin{scriptsize}
\begin{center}
  \caption{Literature overview on works that include applications featuring decision dependent information discovery.}
          \sffamily
\renewcommand{\arraystretch}{1.8}
    \begin{tabular}{|llll|}
        \hline
    \hline
    Citations & Methodology  & Application & Discovery Decisions\\
    \hline
\cite{colvin2008stochastic,colvin2010modeling} &Stochastic programming   & \makecell[l]{Pharmaceutical \\ clinical trial planning}  & Clinical trials \\
\hline
\multirow{2}[1]{*}{\cite{jonsbraaten1998class}}  &   \multirow{2}[1]{*}{Stochastic programming} &\makecell[l]{Offshore gas field \\developments} & Drilling process\\
\cline{3-4}
&  & Production planning  & Production start times\\
\hline
\cite{goel2004stochastic} &Stochastic programming &\makecell[l]{Offshore gas field \\ developments}  &Drilling process\\
\hline
\cite{vayanos2011decision}  & \makecell[l]{Robust optimization,\\Stochastic programming}&\makecell[l]{Offshore gas field \\ developments}   &Drilling process\\
\hline
\cite{solak2010optimization}  &Stochastic programming  &\makecell[l]{Offshore gas field \\ developments}&Drilling process\\
\hline
\cite{fragniere2010operations} &Stochastic programming &\makecell[l]{Workforce \\ capacity planning}  &\makecell[l]{Projects start times,\\ resources allocation}\\
\hline
\cite{macias2020endogenous} &Stochastic programming& \makecell[l]{Humanitarian relief\\ distribution}  & \makecell[l]{UAVs to asses\\ the state of the network} \\
\hline
\cite{vayanos2020robust} &Robust optimization & \makecell[l]{Active preference \\ elicitation}   &Questions for users \\
     \Xhline{1pt}
    \end{tabular}\label{tab:literature}%
    \end{center}
    \end{scriptsize}
    \end{table}


Decision dependent information discovery problems are inherently difficult to solve. This complexity stems primarily from the complexity of the recourse policy. The most successful approximation in the literature is based on $K$-adaptability (\cite{vayanos2020robust}) which can provide good quality solutions if $K$ is sufficiently large, however, at the expense of significant computational time. To obtain good and meaningful solutions when the recourse problem is complex, as it is the case for the orienteering and other logistics problem, exact methods are typically useful to accurately capture the recourse policy. However, approximation schemes can be the only practical solution method for large scale problems. It is therefore also important to improve upon the scalability of the $K$-adaptability. The contribution of this paper can be summarized as follows:
\begin{enumerate}
    \item We propose an exact solution scheme for two-stage robust decision dependent information discovery problems with objective uncertainty. The approach is  a nested decomposition scheme that operates in two levels: $(i)$ an outer layer optimizing over the decisions controlling the timing of information discovery in which we leverage Logic Benders decomposition, and $(ii)$ an inner layer optimizing over the recourse decisions and determines their worst-case cost which we solve using a column-and-constraint generation algorithm. We discuss convergence properties of the algorithm and show that the algorithm is practical and at times outperforms the $K$-adaptability both in solution quality and computational time.
    \item We improve upon the $K$-adaptability formulation proposed in \cite{vayanos2020robust} by strengthening the formulation to achieve greater scalability. To this end, we propose two types of valid inequalities. 
    The first leverages on the symmetry of the recourse problem to strengthen the McCormick inequalities used in the linearization of bilinear terms, and the second provides valid inequalities by identified a-priori “optimistic” realizations of the uncertainty parameters. 
    We emphasize that the  proposed formulation can also be used in the decision independent information discovery setting.     We demonstrate the merits of the strengthened formulation numerically both in the decision dependent and decision decision independent information discovery setting.
    

    \item We use the proposed solution methods on a real problem faced at the Alrijne hospital in the Netherlands that aims to optimize the collection of medicine carrying crates from different departments of the hospital. We formulate the problem  as a two-stage robust orienteering problem with decision dependent information discovery. We solve the problem  using both the exact solution scheme as well as  $K$-adaptability approximation and discuss the interplay between degree of information discovery, the level of adaptability and the setup of the orienteering problem to provide guidelines on how the hospital should optimize its crate collection process.
\end{enumerate}
The remainder of the paper is organized as follows. In Section~\ref{sec:prob_descr} we introduce the two-stage robust decision dependent information discovery problem with objective uncertainty and motivate the connection to the orienteering problem. Section~\ref{sec:exa_algo} presents the exact solution scheme, and Section~\ref{sec:k_adapt} presents the strengthened $K$-adaptability formulation. Section~\ref{section:comp} provides numerical evidence of the usefulness of the proposed methods and Section~\ref{sec:case} discusses the case study of Alrijne hospital. The proofs
of all statements can be found in the Electronic Companion to the paper.

\textit{Notation:} Vector and matrices are represented with bold lowercase and uppercase letters, respectively. We denote be $\bm{e}$ the vector of ones and by $\bm 0$ the vector of zeros. For a generic vector $\bm{x}$, $N_{\bm{x}}$ denotes the size of the vector while $\mcl{N}_{\bm{x}}$ denotes its indices. 
Given two vectors $\bm x$, $\bm y$ $\in \rrr^n$, $\bm{x}\circ\bm{y}$ denotes the Hadamard product. For a logical expression $\mcl{L}$, we define the indicator function $\mathbb{I}[\mcl{L}]$ that takes value 1 if $\mcl{L}$ is true and 0 otherwise.

\section{Problem Description}\label{sec:prob_descr}
We study a two-stage robust optimization problem with \emph{decision depended information discovery with objective uncertainty}. In the classical setting of two-stage robust optimization problems without decision dependent information discovery, after taking the first-stage (or here-and-now) decisions, the decision-maker is able to observe the whole of the uncertain vector $\bm \xi\in\Xi\subseteq \rrr^{N_{\uncV}}$ and utilize this information to decide the recourse (or second-stage) decisions. In contrast, in the decision dependent information discovery setting the first-stage decisions $\bm \discV \in \discVdom \subseteq \bbb^{N_{\discV}}$ control the timing of information discovery. In other words, the recourse decisions $\bm y \in \ySet \subseteq \bbb^{N_{\yvk{}}} $ are able to only adapt to components of the uncertain vector $\bm \xi$ if the decision-maker chooses to observe these components using $\discV$, i.e., if $\discV_i$ = 1, the recourse decisions can adapt to $\bm\xi_i$, however, if $\discV_i$ = 0 the decisions  cannot assimilate $\bm\xi_i$ which is treated as if it materializes after decisions $\yvk{}$ have been decided.  Since each  uncertain parameter is associated with a specific component of $\discV$, we set $N_{\discV} =N_{\unckV{}}$.  The problem is formulated as follows:
\begin{subequations}\label{original_pronblem_full}
\begin{equation}\label{original_prob}
\min_{\bm{w} \in \discVdom }\; \Phi(\discV), 
\end{equation}
with 
\begin{equation}\label{exact_phi_definition}
\Phi\left(\discV\right) = \max_{\bm{\ov{\xi} \in \Xi}}\;\min_{\bm{y} \in \mcl{Y}}\;\max_{\uncV{} \in \uncSet(\discV,\ov{\uncV})}\;\uncV^{\top}\bm{C}\bm{w}+\unckV{}^{\top}\bm{P}\bm{y},
\end{equation}
\end{subequations}
where we assume $\uncSet := \{\uncV \in \rrr^{N_{\xi}} :\bm{A}\uncV \leq \bm{b}\}$ for some matrices $\bm{A} \in \rrr^{L \times N_{\uncV}}$, $\bm{C} \in \rrr^{N_{\uncV}\times N_{\discV}}$, $\bm{P} \in \rrr^{N_{\uncV} \times N_{\yvk{}}}$ and vector $\bm{b} \in \rrr^{L}$.
Problem (\ref{original_pronblem_full}) imposes decision dependent information discovery by utilizing an auxiliary uncertain vector $ \ov{\bm\xi}\in\Xi$ and set $\Xi(\discV,\ov{\bm{\xi}}) := \{\uncV \in \uncSet : \discV \circ \uncV = \discV \circ \ov\uncV \}$ which is a subset of the uncertainty set $\Xi$. Vector $\ov \uncV$ can be thought of as an independent uncertain vector that materializes before decisions $\yvk{}$ are decided, yet, since constraint $\discV \circ \uncV = \discV \circ \ov\uncV$ imposes that $\ov{\uncV}_i = \uncV_i$ if $\discV_i = 1$, the recourse decisions cannot adapt to what has not been observed.  The inner most maximization in~(\ref{exact_phi_definition}) accounts for the worst-case cost for the components of $\uncV$ which have not been observed, i.e., $\discV_i = 0$. For each $\discV$, $\Phi(\discV)$ represents the cost of the recourse problem. For further information, motivation and construction of problem (\ref{exact_phi_definition}) we refer to \cite{vayanos2020robust}.

Although the approaches we develop in this paper are applicable to any problem that can take the form \eqref{original_pronblem_full}, we now focus on an instance of the problem in which decisions $\yvk{} \in \ySet$ encode feasible paths of an orienteering problem. 
The orienteering problem is a routing problem  that aims to generate a routing path through a set of nodes such that it maximizes the profit associated with the visited nodes, while satisfying a  given maximum travel time duration (\cite{tsiligirides1984heuristic,vansteenwegen2011orienteering,gunawan2016orienteering}). The problem is defined by: $(i)$ a complete undirected graph $\mcl{G}(\mcl{N},\mcl{E})$ where  $\mcl{N} := \{0, 1, \ldots, |\mcl{N}|\}$  is the set of nodes and $\edgeSet := \left\{(i,j) \in \nodeSet \times \nodeSet : i < j  \right\}$ is the set of edges, $(ii)$ travel times $\travt{e}{}$  associated with each edge $e = (i,j)\in\mcl{E}$ and a maximum travel time duration $T$, $(iii)$ a start and destination nodes denoted by $\snode$ and $\dnode$, respectively and $(iv)$ a profit $\uncV_i$ for each node $i\in \til{\mcl{N}}:=\mcl{N} \setminus \{s,d\}$. 
The set of all feasible paths can be formally defined as follows
\begin{equation}\label{Y_set}
    \mathcal{Y} = \left\{
    \begin{array}{l@{\,}l}
    \bm y\in\{0,1\}^{|\mathcal{\til N}|}\,:&\,\exists \bm z\in\{0,1\}^{|\mathcal{E}|}\text{ s.t. } \bm t^\top \bm z \leq T,\,\displaystyle\sum_{e \in C(s)}\bm{z}_{e} = \sum_{e \in \mcl{C}(d)}\bm{z}_{e} =1\\
    &\displaystyle \sum_{e \in \mcl{C}(i)}\bm{z}_{e} = 2 \bm{y}_{i}\;\; \forall i \in \mathcal{\widetilde N},\,\displaystyle\sum_{e \in \mathcal{U}(\mcl{S})}\bm{z}_{e} \leq \sum_{i \in \mcl{S}}\bm{y}_{i} - \bm{y}_{u} \;\; u \in \mcl{S},\,\forall \mcl{S} \subseteq \widetilde{N}
    \end{array}
    \right\},
\end{equation}
where, given a subset of nodes $\mcl{I} \subset \nodeSet$ we indicate the edge cutset of set $\mcl{I}$ as $\cutset{\mcl{I}} := \{(i,j) \in \edgeSet : i \in \mcl{I}, j \notin \mcl{I} \lor j \in \mcl{I}, i \notin \mcl{I}  \}$ (i.e., the set of edges $(i,j)$ such that $i \in I,
j \not \in I$ or vice-versa) and with $\mcl{U}(I)$ the set of edges such that the associated two nodes belong to $I$, i.e.,~$\mcl{U}(I):= \{(i,j) \in \edgeSet : i \in \mcl{I}, j \in \mcl{I}\}$. For the sake of notation, we also indicate with $\cutset{i}$ the edge cutset of a
single node $i \in \nodeSet$. In this representation of set $\ySet$, the binary variables $\yvk{}$ indicate which nodes are visited and binary variable $\xvk{}_{e}$ takes value 1 if edge $e$ is traversed and 0 otherwise. 
Although a compact formulation for the orienteering problem is possible by formulating subtours elimination constraints according to the Miller-Tucker-Zemlin formulation (\cite{miller1960integer}), $\ySet$ is defined using a formulation similar to the one used in \cite{fischetti1998solving}. Such a formulation features an exponential number of constraints and therefore, even for the deterministic case must be solved via a branch-and-cut algorithm. However, it typically results in stronger linear relaxations than its compact counterpart.

In the decision dependent information discovery setting, the profits $\uncV_i$ at each node $i \in \til{\mcl{N}}$  are treated as uncertain and are elements of a known uncertainty set $\Xi$. In the first period, the decision maker is able to ``place sensors", corresponding to decisions $\discV$, on nodes so as to observe their associated profit. Using the obtained information, the routing decisions $\yvk{}$ are taken in view of the worst-case profits for those nodes whose profits are not observed. The decision dependent information discovery problem, which we term \emph{sensor placement orienteering problem}, can be written as
\begin{equation}\label{SensorPlacement}
\max_{\discV \in \discVdom}\; \min_{\ov \uncV \in \Xi} \max_{\yvk{} \in \mathcal{Y}}\; \min_{\uncV \in \Xi(\discV,\ov \uncV)}\bm \xi^\top \bm y.
\end{equation}
In the following simple example, we showcase that sensor placement can significantly impact the routing decisions, thus concluding of the need to co-optimizing sensor placements and routing decisions. 

\example \label{example:sol:two_stage_prob} We consider an instance of the orienteering problem where $\mcl{N} = \{0, 1, 2, 3\}$ (0 is both the start and destination node), $\bm{t}_{01} = \bm{t}_{03} = \bm{t}_{02} = 1$, $\bm{t}_{12} = \bm{t}_{13} = \sqrt{2}$, $\bm{t}_{13} = 2$, and the maximum travel time $T = 3.5$. 
In this setting, the feasible routes $\bm{y}$ implied by constraints (\ref{Y_set}) consist of all the paths visiting at most two nodes, i.e., $\mcl{Y} = \{(0,0,0)^\top,(1,1,0)^\top, (0,1,1)^\top, (1,0,0)^\top, (0,1,0)^\top, (0,0,1)^\top\}$.
The profit at each node is represented by the uncertainty vector $\uncV$ that belongs to the uncertainty set $\Xi = \{\bm{\xi} \in \rrr^{3}_+ : \bm{e}^\top\bm{\xi} = 1 \}$, modeling that the total profit (100\%) throughout the network is known, however, its geographical breakdown at each node is uncertain. The sensor placement problem maximizes $\Phi(\discV)$ over the set of admissible sensor placements $\discVdom = \left\{\discV \in \rrr^3 :  \bm{e}^\top\discV \leq B\right\}$, where $B$ is the maximum number of sensors allowed to be placed throughout the network. The profit $\Phi(\discV)$ associated with each placement is 
 \begin{equation}\label{example:original_prob}
\Phi(\discV) = \min_{\bm{\ov{\xi}} \in \Xi}\; \max_{\bm{y} \in \mcl{Y}}\; \min_{\uncV \in \Xi(\bm{w}, \bm{\ov{\xi}})}  \bm{\xi}_{1}\bm{y}_{1} + \bm{\xi}_{2}\bm{y}_{2} + \bm{\xi}_{3}\bm{y}_{3}.
 \end{equation}

For $B = 0$ the only admissible sensor placement is $\discV = \bm{0}$, i.e., no sensors are placed, which implies that $\Xi(\bm w,\bm{\ov \xi}) = \Xi$ for all $\bm{\bar\xi}\in\Xi$. This leads to a robust orienteering problem where the decision maker decides the routing decisions in view of the worst-case realization of $\bm\xi\in\Xi$. As the admissible routes allow to visit at most two nodes, the worst-case allocates all profit to the non-visited node, resulting in the worst-case profit of 0 for any routing decision.

\begin{figure}[htbp!]
\centering
\begin{subfigure}[b]{.45\linewidth}
\centering
\resizebox{.70\linewidth}{!}{
   \begin{tikzpicture}[->,>=stealth',shorten >=2pt, node distance=2.5cm,
                    semithick,framed,background rectangle/.style={thick,fill=gray1,draw=black}] 
    \tikzstyle{every state}=[fill=white,draw=black,text=black,font=\footnotesize]

 \node[triangle, draw] (dep) at (-10,0) {\footnotesize{0}};
      \node[shape=circle,draw,
    left of = dep
    ] (n1) {1};
      \node[shape=circle,draw, above of = dep
      ] (n2) {2};
\node[shape=circle,draw, right of = dep
] (n3) {3};
 \draw[->](dep)  to node [auto] {} (n1); 
  \draw[->]  (n1) to node [auto] {} (n2);
    \draw[->]  (n2) to node [auto] {} (dep);


\end{tikzpicture} 
}
\subcaption{Optimal route if $\uncV_1 + \uncV_2 \geq \uncV_2 + \uncV_3$.}\label{subfig:example1:we:sub1}
\end{subfigure}
\begin{subfigure}[b]{.45\linewidth}
\centering
\resizebox{.70\linewidth}{!}{
    \begin{tikzpicture}[->,>=stealth',shorten >=2pt, node distance=2.5cm,
                    semithick,framed,background rectangle/.style={thick,fill=gray1,draw=black}] 
    \tikzstyle{every state}=[fill=white,draw=black,text=black,font=\footnotesize]

 \node[triangle, draw] (dep) at (0,0) {\footnotesize{0}};
      \node[shape=circle,draw,
    left of = dep
    ] (n1) {1};
      \node[shape=circle,draw, above of = dep
      ] (n2) {2};
\node[shape=circle,draw, right of = dep
] (n3) {3};
 \draw[->](dep)  to node [auto] {} (n2); 
  \draw[->]  (n2) to node [auto] {} (n3);
    \draw[->]  (n3) to node [auto] {} (dep);

\end{tikzpicture} 
}
\subcaption{Optimal route if  $\uncV_1 + \uncV_2 \leq \uncV_2 + \uncV_3$.}\label{subfig:example1:we:sub2}
\end{subfigure}
\caption{Routing policy for $\discV = (1,1,1)^\top$, collecting in the  worst-case 50\% of the total profit.}\label{fig:example1:we}
\end{figure}
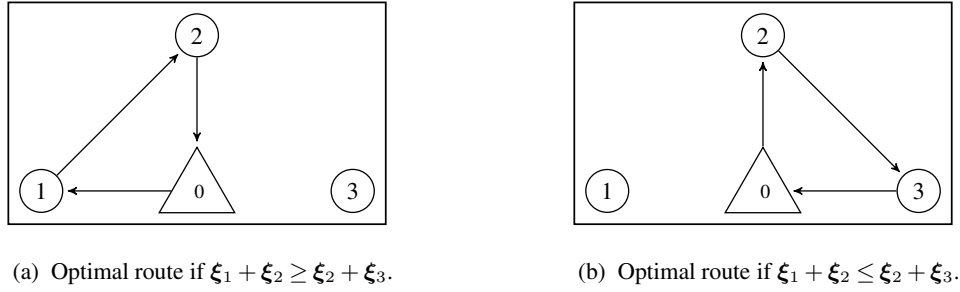

For $B = 3$, $\discV = (1,1,1)^\top$ is admissible which implies placing sensors in all nodes, thus $\Xi(\discV,\ov\uncV) = \{\ov \uncV \}$ for all $\ov \uncV\in\Xi$. Intuitively, $\discV = (1,1,1)^\top$ is an optimal solution, as the routing decisions will be chosen in full knowledge of the worst-case profits in the network,~i.e., more information is better for the recourse decisions. This is formalized in Lemma~\ref{lemma_exa_function_prop} in the following section. As demonstrated in Figure~\ref{fig:example1:we}, the routing is a policy of the realizations of $\bm\xi$.
If $\uncV_1 + \uncV_2 \geq \uncV_2 + \uncV_3$  (Figure~\ref{subfig:example1:we:sub1}), nodes 1 and 2 are visited, while if $\uncV_1 + \uncV_2 \leq \uncV_2 + \uncV_3$ the decision-maker select the path that visits node 2 and node 3  (Figure~\ref{subfig:example1:we:sub1}). 
In this case, the adversarial nature allocates  profit such that $\uncV_1 + \uncV_2 = \uncV_2 + \uncV_3$, resulting in the worst-case profit of 0.5, or 50\% of the total profit.

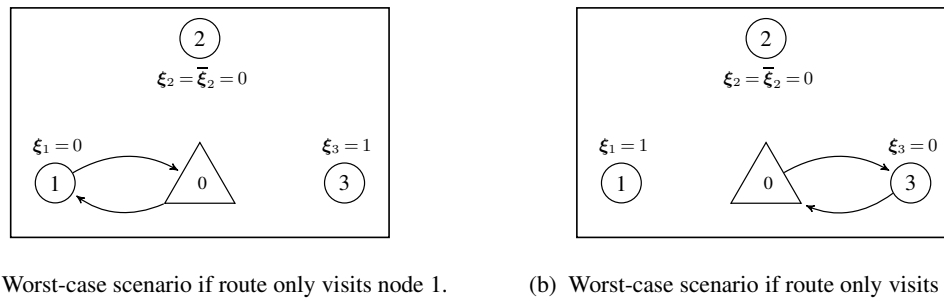
\begin{figure}[htbp!]
  \centering
   \captionsetup[subfigure]{font=footnotesize,labelfont=footnotesize}
\begin{subfigure}[b]{.45\linewidth}
\centering
\resizebox{.70\linewidth}{!}{
    \begin{tikzpicture}[->,>=stealth',shorten >=2pt, node distance=2.5cm,
                    semithick,framed,background rectangle/.style={thick,fill=gray1,draw=black}] 
    \tikzstyle{every state}=[fill=white,draw=black,text=black,font=\footnotesize]

 \node[triangle, draw] (dep) at (0,0) {\footnotesize{0}};
      \node[shape=circle,draw,
    left of = dep,label={[align=left]above:\footnotesize{$\uncV_1 = 0$}} ] (n1) {1};
      \node[shape=circle,draw, above of = dep,label={[align=left]below:\footnotesize{$\uncV_2 = \ov{\uncV}_2= 0$}}] (n2) {2};
\node[shape=circle,draw, right of = dep,label={[align=left]above:\footnotesize{$\uncV_3 = 1$}}] (n3) {3};
 \draw[bend left,->](dep)  to node [auto] {} (n1); 
  \draw[bend left,->]  (n1) to node [auto] {} (dep); 
\end{tikzpicture} 
}
\caption{Worst-case scenario if route only visits node 1.}\label{subfig:example1:w_2:sub1}
\end{subfigure}
\begin{subfigure}[b]{.45\linewidth}
\centering
\resizebox{.70\linewidth}{!}{
    \begin{tikzpicture}[->,>=stealth',shorten >=2pt, node distance=2.5cm,
                    semithick,framed,background rectangle/.style={thick,fill=gray1,draw=black}] 
    \tikzstyle{every state}=[fill=white,draw=black,text=black,font=\footnotesize]
 \node[triangle, draw] (dep) at (0,0) {\footnotesize{0}};
 \node[shape=circle,draw,left of = dep,label={[align=left]above:\footnotesize{$\uncV_1 = 1$}} ] (n1) {1};
\node[shape=circle,draw, above of = dep,label={[align=left]below:\footnotesize{$\uncV_2 = \ov{\uncV}_2= 0$}}] (n2) {2};
\node[shape=circle,draw, right of = dep,label={[align=left]above:\footnotesize{$\uncV_3 = 0$}}] (n3) {3};
\draw[bend left,->](dep)  to node [auto] {} (n3); 
\draw[bend left,->]  (n3) to node [auto] {} (dep); 
\end{tikzpicture} 
}
\caption{Worst-case scenario if route only visits node 3.}\label{subfig:example1:w_2:sub2}
\end{subfigure}
\caption{Worst-case scenarios for $\discV = (0,1,0)^\top$,~i.e., observing $\uncV_2$. All routing policies result to a worst-case profit of 0.}\label{fig:example1:w_2}
\end{figure}

\begin{figure}[htbp!]
  \centering
    \captionsetup[subfigure]{font=footnotesize,labelfont=footnotesize}
\begin{subfigure}[b]{.45\linewidth}
\centering
\resizebox{!}{.40\linewidth}{
    \begin{tikzpicture}[->,>=stealth',shorten >=2pt, node distance=2.2cm,
                    semithick,framed,background rectangle/.style={thick,fill=gray1,draw=black}] 
    \tikzstyle{every state}=[fill=white,draw=black,text=black,font=\footnotesize]

 \node[triangle, draw] (dep) at (0,0) {\footnotesize{0}};
      \node[shape=circle,draw,
    left of = dep,
   label={[align=left]above:\footnotesize$\ov{\uncV}_1 = \uncV_1$} 
    ] (n1) {1};
      \node[shape=circle,draw, above of = dep
      ] (n2) {2};
\node[shape=circle,draw, right of = dep
] (n3) {3};
 \draw[bend left,->](dep)  to node [auto] {} (n1); 
  \draw[bend left,->]  (n1) to node [auto] {} (dep); 
\end{tikzpicture} 
}
\caption{Optimal route if  $\uncV_{1} \geq 0.5 
$.}\label{subfig:example1:w_1:sub1}
\end{subfigure}
\begin{subfigure}[b]{.45\linewidth}
\centering
\resizebox{!}{.40\linewidth}{
    \begin{tikzpicture}[->,>=stealth',shorten >=2pt, node distance=2.5cm,
                    semithick,framed,background rectangle/.style={thick,fill=gray1,draw=black}] 
    \tikzstyle{every state}=[fill=white,draw=black,text=black,font=\footnotesize]

 \node[triangle, draw] (dep) at (0,0) {\footnotesize{0}};
      \node[shape=circle,draw,
    left of = dep
   ,label={[align=left]above:\footnotesize$\ov{\uncV}_1 = \uncV_1$} 
    ] (n1) {1};
      \node[shape=circle,draw, above of = dep
      ] (n2) {2};
\node[shape=circle,draw, right of = dep
] (n3) {3};
 \draw[->](dep)  to node [auto] {} (n2); 
  \draw[->](n2)  to node [auto] {} (n3); 
  \draw[->] (n3) to node [auto] {} (dep); 
\end{tikzpicture} 
}
\caption{Optimal route if  $\uncV_{1} \leq 0.5$
.}\label{subfig:example1:w_1:sub2}
\end{subfigure}
\caption{Routing policy for $\discV = (1,0,0)$, i.e., observing $\uncV_1$, collecting in the  worst-case 50\% of the total profit.}\label{fig:example1:w_1}
\end{figure}
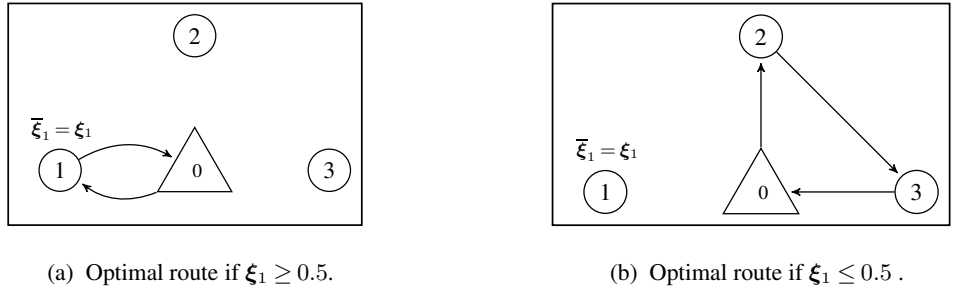

  
 When $B = 1$, there are four admissible sensor placements $\bm w = (0,0,0)^\top$, $\bm w = (1,0,0)^\top$, $\bm w = (0,1,0)^\top$ and $\bm w= (0,0,1)^\top$. For $\bm w = (0,1,0)^\top$,~i.e.,~placing  a sensor only in node 2, regardless of the optimal routing policy the worst-case profit is  0. This is the case as the worst-case profit realization of node 2 is zero, which is represented by $\ov{\uncV}_2 = 0$, implying that $\Xi(\discV,\ov\uncV) = \{\uncV\in\mathbb{R}^3_+: \uncV_1 + \uncV_3 \leq 1,\; \uncV_2 = \ov{\uncV}_2= 0\}$. Thus, the total profit can be placed either in node 1 or node 3, however, since the routing decisions are lacking this information and since the construction of $\mathcal{Y}$ allows only one of nodes 1 and 3 to be visited, in the worst-case none of the profit is collected as shown in  Figure~\ref{fig:example1:w_2}. We now consider placing a sensor on node 1, $\bm w = (1,0,0)^\top$. Due to the maximum time duration $T = 3.5$, the admissible routing decisions are allowed to visit either nodes 1, or nodes 2 and 3. The optimal routing policy is: if $\uncV_1 \geq 0.5$ then visit node 1, while if $\uncV_1 \leq 0.5$ then visit nodes 2 and 3, as shown in Figure~\ref{fig:example1:w_1}. The corresponding worst-case realization sets $\ov \uncV_1 = 0.5$, implying $\Xi(\discV, \ov{\uncV}) = \{\uncV\in\mathbb{R}^3_+ : \uncV_2 + \uncV_3 = 0.5,\, \uncV_1 = \ov \uncV_1 =0.5\}$, which results in a worst-case profit of 0.5. Note that due to the topology of the network, the solutions for $\bm w = (1,0,0)^\top$ and $\bm w  = (0,0,1)^\top$ are symmetric thus they achieve the same objective value. Since $B = 1$ and $B = 3$ both achieve an optimal profit of 0.5, $B = 2$ will  result in the same worst-case profit.
 
 
\endexample
The example demonstrates that optimally placing sensors on only a subset of the nodes is sometimes sufficient for achieving the same optimal worst-case profit as placing sensors on all nodes. However, misplacing sensors can lead to the same worst-case profit as having no sensors in the network. Therefore, there is a need for optimizing over the sensor placement which can have a significant impact on the routing policies. In the following section we propose an exact procedure for optimally solving problem~\eqref{original_pronblem_full}.

\section{Exact Algorithm for Problem (\ref{original_pronblem_full})}\label{sec:exa_algo}
In this section we propose an exact algorithm for solving problem~\eqref{original_pronblem_full}. Utilizing ideas from Logic Benders we iteratively underestimate $\Phi(\discV)$ using linear inequalities.
To generate an inequality, the algorithm requires to exactly evaluate $\Phi(\bm w)$  for each $\discV$, and since \eqref{exact_phi_definition} is in essence a two-stage linear robust optimization we design a column-and-constraint generation algorithm to do so. In Section~\ref{sec:LogicBenders} we discuss the properties of $\Phi(\discV)$ and propose valid inequalities and in Section \ref{sec:exa:sub_problem} we present the column-and-constraint generation algorithm for evaluating $\Phi(\discV)$. Section~\ref{sec:exa_extension} in the Electronic Companion discusses extensions to accommodate a richer class of two-stage robust optimization problems with decision dependent information discovery. 

\subsection{Logic Benders Decomposition for Problem (\ref{original_prob})}\label{sec:LogicBenders}
Function $\Phi(.)$ is a real valued function which is in general not convex on $[0, 1]^{N_{\discV}}$. We are however interested in the values of the function at the binary points $\bbb^{N_{\discV}}$. Our decomposition leverages on the fact that the convex lower envelope of $\Phi(.)$ is a convex piece-wise linear function on $[0, 1]^{N_{\discV}}$ which exactly approximate $\Phi(.)$ at all the binary points \cite[Theorem 1]{zou2019stochastic}.
The convex lower envelope can be described by the following set of inequalities
\begin{subequations}\label{exact_convex_envelope}
\begin{equation}\label{exact_integer_cut}
\Phi(\discV) \geq \Phi(\discV') - \left(\Phi(\discV') - \un{\Phi}\right)\phi(\discV', \discV) \quad \forall \discV'\; \in \discVdom,
\end{equation}
where
\begin{equation}\label{exact_function}
\phi(\discV', \discV) = \sum_{\substack{i\; \in\; \mcl{N}_{\discV} :\; \discV'_i = 1}}(1 - \discV_i) +  \sum_{\substack{i\; \in\; \mcl{N}_{\discV} :\; \discV'_i = 0}}\discV_i,
\end{equation}
and $\un{\Phi} \in \rrr$ such that $\un{\Phi} \leq \Phi(\discV)$ for all $\discV \in \discVdom$, i.e., any lower bound on $\Phi(\discV)$.
\end{subequations}
This family of inequalities is one of the main ingredient of Logic Benders based decomposition algorithms. Logic Benders decomposition was formally introduced in \cite{hooker2003logic}   and is usually applied to two-stage, e.g.,  \cite{caroe1998shaped}, \cite{laporte1993integer}, and multi-stage   \cite{zou2019stochastic} stochastic optimization problems. 
However, they are not widely adopted within robust optimization where decomposition approaches based on column-and-constraints generation are more common, e.g.,  \cite{zeng2013solving,ayoub2016decomposition}.


{\proposition{Inequalities \eqref{exact_convex_envelope} are valid and tight at each point  $\discV\; \in \discVdom$.}\label{preposition_tight_and_valid}
\endproposition}

In line with Benders-like decomposition scheme, we leverage on Proposition \ref{preposition_tight_and_valid} to define the following relaxation of problem \eqref{original_prob}.
\begin{equation}\label{exact_master}
    \Theta_{\widehat{\discVdom}} =\left[ \begin{array}{lll}
     \min & \varphi &\\
      \st& \varphi \in \rrr,\ \discV \in \mcl{W}&\\
        & \varphi \geq \Phi(\discV') - \left(\Phi(\discV') - \un{\Phi}\right)\phi(\discV', \discV) & \quad\forall\; \discV' \in \widehat{\discVdom}
    \end{array}\right],
\end{equation}
where $\widehat{\discVdom}\subseteq \discVdom$. Problem~\eqref{exact_master} includes only a subset of all the possible 
inequalities~\eqref{exact_convex_envelope}, thus $\Theta_{\widehat{ 
\discVdom}} \leq \Phi(\discV)$ for all $\widehat{\discVdom}\subseteq \discVdom$. 
Algorithm~\ref{Algo_Logic_Benders} optimally solves problem~\eqref{original_prob} by identifying 
which constraints~\eqref{exact_convex_envelope} the progressive solution $\discV^*$ violates and iteratively add them to set $\widehat{\discVdom}$. 
\begin{algorithm}[H]\label{Algo_Logic_Benders}
\SetAlgoLined
\SetKwInOut{Input}{Input}
\SetKwInOut{Output}{Output}
\Output{Optimal solution $\discV^*$ and optimal solution value $\Phi(\discV^*)$.}
\textbf{Initialize:} $ ub  =+\infty, \, lb  =-\infty$, $\widehat{\discVdom} = \emptyset$;\\
 Solve Problem (\ref{exact_master}) and let $\discV^* $ be the optimal solution. Set $lb$ = $           \Theta_{\widehat{\discVdom}}$\;\label{StepSolveMP} 
 Evaluate $\Phi(\discV^*)$ by solving~\eqref{exact_phi_definition} and set  $ub$ = $\Phi(\discV^*)$\;\label{StepEvaluatePhi} 
 \lIf{ub  -  lb   $>$ 0}
 {
 $ \widehat{\discVdom} = \widehat{\discVdom}\; \cup\; \{\discV^{*}\} $ and go to Step~\ref{StepSolveMP}
  }
 \lElse
 {
 \Return $\discV^{*}$ and  $\Theta_{\widehat{\discVdom}}$
 }
\caption{Exact Algorithm for Problem \eqref{original_prob}}
\end{algorithm}\noindent
At each iteration, problem~\eqref{exact_master} is 
solved to identify the progressive solution $\discV^*$, and the value of the recourse problem $\Phi(\discV^*)$ is evaluated by solving problem \eqref{exact_phi_definition}. The algorithm terminates when the optimal value  $\Theta_{\widehat{\mathcal{W}}}$ of the relaxed problem and $\Phi(\discV^*)$ coincide, otherwise $\discV^*$ is added to set $\widehat{\discVdom}$ and the process is repeated. The following theorem demonstrates that the algorithm terminates in finite number of steps.

\begin{theorem}\label{th_exa_converge}
Algorithm \ref{Algo_Logic_Benders} converges to the optimal solution of problem (\ref{original_prob}) in a finite number of iterations.
\end{theorem}
In line with modern implementations of Benders-like algorithms (see \cite{adulyasak2015benders,fischetti2017redesigning}), rather than executing  Algorithm (\ref{Algo_Logic_Benders}), we embed the resolution of the relaxed master problem  (\ref{exact_master}) in a branch-and-cut algorithm. Within this algorithm, every time a feasible solution  $(\varphi^*,\ \discV^*)$ is encountered in the search tree, function  $\Phi(\discV^*)$ is evaluated and if $\varphi^* < \Phi(\discV^*)$, the corresponding cut is added to problem (\ref{exact_master}). The branch-and-cut algorithm ends when all the nodes are pruned and an optimal solution to problem (\ref{exact_master}) is achieved. In other words, with respect to Algorithm \ref{Algo_Logic_Benders}, we execute Step \ref{StepSolveMP} only one time while Step \ref{StepEvaluatePhi} will be execute every time a new feasible solution $\discV \in \discVdom$ is founded in the search tree.  

Inequalities~\eqref{exact_convex_envelope}  are valid and tight at every $\discV \in \discVdom$, therefore they ensure correctness and convergence of the algorithm. However, they can be, in general, extremely weak at fractional solution, i.e., $\discV \in [0,1]^{N_{\discV}}$ slowing down the resolution of problem \eqref{exact_master}. This is a typical draw-back of Logic Benders solution schemes (\cite{hooker2003logic,zou2019stochastic}). In the following, we leverage on the fact that $\discV$ encodes the amount of information available to the recourse decisions thus we are able to better describe the behavior of $\Phi(\discV)$. This will allow us to build tighter inequalities , which we term \emph{information inequalities}, achieving faster convergence rate in practice. 
\subsubsection{Information Inequalities}
We now consider a variant of problem \eqref{original_pronblem_full} where the costs affecting discovery variables are deterministic rather than uncertainty, i.e., $\uncV^{\top}\bm{C}\bm{w} = \bm{c}^\top\discV$ for a given vector $\bm{c} \in \rrr^{N_{\bm{w}}}$. In this setting, problem \eqref{original_prob} can be written as $\displaystyle\min_{\bm w\in\mathcal{W}}\bm{c}^\top\discV + \Phi(\bm w)$, where  $\displaystyle\Phi(\bm w) = \max_{\bm{\bar\xi} \in \uncSet} \min_{\bm y\in\mathcal{Y}} \max_{\bm \xi\in\Xi(\bm w, \bm{\bar \xi})} \bm \xi^\top \bm P\bm y$. The following lemma captures the intuition that as the information increases,~i.e., more entries of $\discV$ take values equal to 1, the decision-maker can take more informed recourse decisions, thus the value of the $\Phi(\bm w)$ is non-increasing.
{
\lemma{}\label{lemma_exa_function_prop}
Let $\discV^{1}$, $\discV^2$ $\in \discVdom$. \mbox{If $\discV^1 \geq \discV^2$ (where the inequality holds component wise), then $\Phi(\discV^1) \leq \Phi(\discV^2)$.}
\endlemma}
It is important to note that Lemma~\ref{lemma_exa_function_prop} does not hold when the cost associated with $\discV$ is uncertain. In this setting, the value of $\Phi(\discV)$ can increase as more components of $\discV$ take the value one. Intuitively, this can occur as the worst-case $\uncV$ is chosen to maximize $\uncV^\top \bm{C}\discV+ \uncV^\top \bm{C}\yvk{}$, hence the it can shift more weight on $\uncV^\top \bm{C}\discV$. 

Leveraging on Lemma \ref{lemma_exa_function_prop}, we propose the following \emph{information inequalities}
\begin{subequations}\label{exact_information_cut}
\begin{equation}\label{exact_lift_integer_cut}
\Phi(\discV) \geq \Phi(\discV') - \left(\Phi(\discV') - \un{\Phi}\right)\rho(\discV', \discV) \quad \forall \discV' \in \discVdom,
\end{equation}
where
\begin{equation}\label{exact_lift_function}
\rho(\discV', \discV) = \sum_{\substack{i \in \mcl{N}_{w} : \discV'_i = 0}}\discV_i.
\end{equation}
\end{subequations}
The following proposition proves the validity of the inequalities.
{\proposition
Inequalities~\eqref{exact_information_cut} are valid and tight at every feasible integer solution $\discV \in \discVdom$ and at least as tight as inequalities~\eqref{exact_convex_envelope} at every fractional solution $\discV'' \in [0,1]^{N_{\discV}}$. \label{proposition_information_inequalities}
\endproposition}

As we need to evaluate $\Phi(\discV)$ exactly at each iteration of the algorithm, solving efficiently problem~\eqref{exact_phi_definition} is crucial for a practical implementation. In the next section we propose a column-and-constraint generation algorithm for evaluating $\Phi(\discV)$.


\subsection{Column-and-Constraint Generation Algorithm for Problem (\ref{exact_phi_definition}) }\label{sec:exa:sub_problem}
Problem (\ref{exact_phi_definition}) is a two stage robust optimization problem with objective uncertainty, where the role of the first and second stage decisions are taken by $\ov \uncV$ and $\uncV$, respectively, and the role of the ``uncertainty" is taken by the vector $\yvk{}$. 
To solve this problem, we develop a column-and-constraint generation algorithm. We express $\displaystyle\Phi(\discV) = \max_{\ov \uncV \in \uncSet} \psi(\discV,\ov \uncV)$, where  
\begin{equation}\label{subproblem_ccg}
    \psi(\discV, \ov\uncV) = \min_{\yvk{} \in \ySet}\;\max_{\uncV \in \uncSet(\discV, \ov\uncV)}\uncV^{\top}\bm{C}\bm{w}+\unckV{}^{\top}\bm{P}\bm{y}.
\end{equation}
The following problem constitutes a relaxation of problem~(\ref{exact_phi_definition})
\begin{equation}\label{master_ccg}
\begin{array}{lll}
 \Psi_{\widehat\ySet}(\discV) &= &\displaystyle \max_{\ov\uncV \in \uncSet}\; \min_{\yvk{} \in \widehat{\ySet}}\;\max_{\uncV \in \uncSet(\discV, \bar\uncV)}\bm\xi^\top \bm{C}\bm{w}+\unckV{}^{\top}\bm{P}\bm{y}\\[3ex]
 &= &
 \left[\begin{array}{rll}
 \max& \tau  \\
 \st&\tau \in \rrr,\ \ov{\uncV{}} \in \uncSet,\ \uncV:\ySet\mapsto\rrr^{N_{\uncV}}&\\
 & \tau \leq \uncV(\bm y)^{\top}\bm{C}\bm{w}+\bm \xi(\bm y)^{\top}\bm{P}\bm{y} & \forall \yvk{} \in \widehat\ySet \\
 &\unckV{}(\yvk{}) \in \Xi(\bm{w}, \bm{\ov{\xi}}) & \forall \yvk{} \in \widehat\ySet
\end{array} \right], 
\end{array}
\end{equation}
where  $\widehat \ySet \subseteq \ySet$, and the second equality follows from \cite[Theorem 1]{vayanos2020robust}. For any subset $\widehat \ySet $, the optimal value $\Psi_{\widehat \ySet}(\bm w)$ is an upper bound on $\Phi(\discV)$. 

Algorithm~\ref{Algo_CCG} solves problem~\eqref{exact_phi_definition} to optimality. The algorithm operates by iteratively solving the relaxed problem~\eqref{master_ccg} to obtain the upper bound $\Psi_{\widehat{\ySet}}(\discV)$ and get the progressive solution $\ov \uncV$, whose value is subsequently evaluated by solving problem~\eqref{subproblem_ccg} to get the lower bound $\psi(\discV,\ov{\uncV})$ and the value $\yvk{*}$ that achieves it. If the bounds coincide, then the algorithm terminates, otherwise $\yvk{*}$ is added in set $\widehat \ySet$ and the process is repeated. Notice that for given $\discV$ and $\ov \uncV$ problem \eqref{subproblem_ccg} is a static robust optimization problem and can be either reformulated as a single mixed-integer linear program by dualizing the inner maximization, or can be solve through a constraint generation algorithm, see \cite{bertsimas2016reformulation} for a comparison between the two approaches. In our work, we follow the former approach. The algorithm is termed column-and-constraint generation as for every $\yvk{*}$ added in set $\widehat \ySet$ new constraints are added in problem~\eqref{master_ccg}, as well as new decision variables (``columns") through $\uncV(\yvk{})$. 
The following theorem demonstrates that the algorithm terminates in finite number of steps.


\begin{algorithm}[H]\label{Algo_CCG}
\SetAlgoLined
\SetKwInOut{Input}{Input}
\SetKwInOut{Output}{Output}
\Input{A vector $\discV \in \discVdom$.}
\Output{Value of $\Phi(\discV)$.}
\textbf{Initialize:} $ lb  =-\infty, \, ub  =\infty$, $\mathcal{\widehat Y}= \emptyset$;\\
 Solve (\ref{master_ccg}) and let $\bm{\ov{\xi}^*}$ be the optimal solution. Set $ub$ = $\Psi_{\widehat\ySet}(\discV)$\;\label{StepSolveSP} 
Evaluate $\psi(\bm{w},\ \bm{\ov{\xi}^*})$ by solving problem (\ref{subproblem_ccg}) and let $\yvk{*}$ the optimal solution. Set $lb$ =  $\psi(\bm{w},\ \bm{\ov{\xi}^*})$\;
 \lIf{ub  -  lb   $>$ 0}
 {
 $ \widehat{\ySet} = \widehat{\ySet}\; \cup\; \{\yvk{*}\} $ and go to Step \ref{StepSolveSP}
  }
 \lElse
 {
 \Return $\Psi_{\widehat\ySet}(\discV)$
 }
\caption{Exact Column-and-Constraint Generation Algorithm for Problem (\ref{exact_phi_definition})}
\end{algorithm}

 \begin{theorem}\label{th_ccg_convergence}
 Algorithm \ref{Algo_CCG} converges after a finite number of iterations. 
 \end{theorem}

\section{Improving the Performance of the \texorpdfstring{$K$}{}-adaptability Approximation}\label{sec:k_adapt}
The two-stage robust optimization problem with decision-dependent information discovery is known to be NP-hard, thus exact algorithms are expected to struggle as the size of the problem instance increases.  To gain tractability an alternative approach is to approximate problem~\eqref{original_pronblem_full}. First introduced by \cite{bertsimas2010finite}, $K$-adaptability is a popular approximation for the class of problems with binary recourse decisions. The approximation reduces the complexity by approximating the recourse decisions $\yvk{}$ with $K$ constant decisions and treating them as first-stage (here-and-now) decisions. After uncertainty materializes (or partly materializes in the decision-dependent discovery setting) the best amongst the $K$ constant decisions is used so as to minimize the worst-case cost. For $\discV = \bm{e}$, i.e., the decision-independent information discovery, \cite{hanasusanto2015k} and \cite{subramanyam2019k} provide efficient solution techniques based on a monolith mixed integer
linear formulation and an iterative solution scheme involving disjunctive programming, respectively. For the decision dependent information discovery case, \cite{vayanos2020robust} also provides a monolithic formulation
which is similar in spirit to the formulation provided in \cite{hanasusanto2015k}. In this section,
we improve upon the formulation proposed in \cite{vayanos2020robust} by strengthening the formulation thus
providing significantly improvements in terms of scalability and computational time. Our results can also be applied in the decision-independent case which is subsumed by the decision depended case. In the following, we first briefly present the $K$-adaptability formulation and then we discuss additional valid inequalities that strengthen the formulation.

The $K$-adaptability approximation of problem (\ref{original_prob}), as introduced in \cite{vayanos2020robust}, is stated as follows:
\begin{equation*}\label{approx_prob}
\begin{array}{>{\displaystyle}l>{\displaystyle}l} 
\min& \max_{\bm{\ov{\xi}} \in \Xi}\; \min_{k \in \mcl{K}}\;\left\{ \max_{\xi \in \Xi(\bm{w}, \bm{\ov{\xi}})} \uncV^{\top}\bm{C}\bm{w} + \bm{\xi}^{\top}\bm{P}\bm{y}^k\right\} \\
\displaystyle\st&\bm{w} \in \mcl{W}, \bm{y}^{k} \in \mcl{Y},\; \forall k \in \mcl{K},
\end{array}
\end{equation*}
where $\mcl{K} = \{1, 2, \ldots, K\}$. The here-and-now variables $\yvk{k}, k \in \Kset$ encode the $K$ constant decisions that define the recourse policy. In view of vector $\ov{\uncV}$ the inner minimization calls for selecting the best piece among the $K$ constant decisions. The higher the value of $K$, the better the $K$-adaptability approximation, while setting $K = |\ySet|$ ensures the optimal solution of \eqref{original_pronblem_full} is recovered. In practice, it turns out that for some problem relatively small values of $K$, can provide good quality solutions. If we contrast $K$-adaptability and the exact algorithm presented in Section~\ref{sec:exa_algo}, $K$-adaptability only uses $K$ constant decisions while Algorithm~\ref{Algo_CCG} dynamically increases the number of constant decisions in $\widehat{\ySet}$ until a set that guarantees optimality is achieved. The $K$-adaptability counterpart can be reformulated as the following mixed integer bilinear programming problem
\cite[Theorem~2]{vayanos2020robust} 
\begin{subequations}\label{mod_bi_problem}
\begin{alignat}{2}
\min\quad&  \bm{b}^\top\left(\bm{\beta} + \sum_{k \in \mcl{K}}\bm{\beta}^{k}\right) 
& & \\
\st \quad &  \bm{w} \in \mcl{W},\; \bm{y}^{k} \in \mcl{Y}\;,  k \in \Kset &&\\
&\bm{\alpha} \in \rrr^{K}_{+},\; \bm{\beta} \in \rrr^{L}_{+},\;\bm{\beta}^{k} \in \rrr^{L}_{+},\;  \bm{\gamma}^{k} \in \rrr^{N_{w}},\; k \in \Kset  & &\\
 & \bm{e}^{\top}\alDV = 1  & &\label{mod_bil:cnstr:sum_alpha}\\
& \bm{A}^\top\bm{\beta} = \sum_{k \in \mcl{K}}\bm{\gamma}^{k}\circ\bm{w}  & &\label{mod_bil:cnstr:balance}\\
&  \bm{A}^\top\bm{\beta}^{k} +\bm{\gamma}^{k}\circ\bm{w} =  \alDV_{k}\left(\bm{C}\discV +\bm{P}{\bm{y}^{k}}\right)  &\quad\forall k \in \mcl{K}.&\label{mod_bil:cnstr:balance_k}
\end{alignat}
\end{subequations}

As shown in \cite{vayanos2020robust}, the bilinear terms $\bm{\gamma}^{k}\circ\bm{w}$, $\bm{\alpha}_{k}{\bm{y}^{k}}$ and $\bm{\alpha}_{k}\discV$, for each $k \in \Kset$, can be linearized using the well known McCormick inequalities (\cite{mccormick1983nonlinear}). Unfortunately, the resulting linearization   is in general characterized by a weak linear relaxation. On top of that, $K$-adaptability reformulation presents large amount of symmetry since the $K$ candidate policies can be permuted to achieve equivalent solutions. Consequently, the $K$-adaptability problem is already challenging to solve for $K \in \{3, 4\}$ (\cite{hanasusanto2015k}). To overcome these issues, we propose several valid inequalities to both strengthen the formulation and decrease the degree of symmetry. In Section~\ref{subsec:simmetry} we propose a new family of inequalities to reduce the symmetry of the problem and in Section \ref{sssec:stronger_lin}, we show how we can strengthen the McCormick inequalities by restricting the domain of variables involved in the linearization. Furthermore, in Section~\ref{sss:lb_cut}, we propose the \emph{optimistic inequalities} defined on optimistic realizations of the uncertainty vector $\uncV$. Finally, in Section~\ref{sss:rlt_cuts} we discuss additional valid inequalities that can be added by exploiting the Reformulation-Linearization-Technique.

\subsection{Addressing Symmetry}\label{subsec:simmetry}
A drawback of the $K$-adaptability formulation is the symmetry with respect to the policies $\yvk{k}$. Indeed, in problem (\ref{mod_bi_problem}) the index of  variables $(\bm \alpha_k,\bm \beta^k,\bm \gamma^k,\bm y^k)$ can be permuted without changing the optimal value of the problem.
To reduce symmetry in the $K$-adaptability reformulation, we need to ensure that permutating values of variables associated to the $K$ policies, does not lead to equivalent solutions. We enforce this condition by imposing that $\alDV$ belongs to set $\mcl{A}_S$ defined as follows
\begin{equation}\label{symmetry_breaking}
\mcl{A}_S  := \left\{ \alDV \in \rrr^{K}_+: \bm{e}^\top\alDV = 1,\;  \bm{\alpha}_k \geq \bm{\alpha}_{k+1}\; \forall k \in \mcl{K} \setminus \{K\} \right\}.
\end{equation}
The following theorem ensures that the addition of constraints defining $\mcl{A}_S$ does not result in a conservative approximation. 
\begin{theorem}\label{th_symmetry_breaking}
At least one optimal solution of problem (\ref{mod_bi_problem}) satisfies  $\alDV \in \mcl{A}_S$.  
\end{theorem}

 In \cite{vayanos2020robust} a different family of symmetry breaking constraints is introduced  which breaks  symmetry by imposing that $\bm y^k,\; k \in \Kset$ satisfy a lexicographic decreasing order. In the next subsection, we will use the properties of set $\mcl{A}_S$ to further strengthen the McCormick inequalities used to linearize the bilinear terms in problem (\ref{mod_bi_problem}). Unfortunately, the symmetry breaking constraints from \cite{vayanos2020robust} cannot be used together with the strengthened McCormick inequalities.

\subsection{Strengthening McCormick Inequalities}\label{sssec:stronger_lin}
A common approach to linearize products of binary and continuous variables is through the so-called McCormik inequalities (\cite{mccormick1983nonlinear}). For every product of two original variables, a new continuous variable and four constraints are introduced. Consider bilinearity $\bm{\alpha}_{k}{\bm{y}^{k}}$ in problem~\eqref{mod_bi_problem}. If $ \ell_{\bm{\alpha}_k}$ and $u_{\bm{\alpha}_k}$ are the lower and upper bounds for variables  $\alDV_k$ for each $k \in \Kset$, respectively, and if we set $\tyvk{k}_i =\alDV_k\yvk{k}_i$, for each $i \in \mathcal{N}_{\bm{y}}$, the corresponding McCormick inequalities state as follows
\begin{equation}\label{McCormik_explicit}
\begin{aligned}
 &\bm{\til{y}}^k_i  \geq \ell_{\bm{\alpha}_{k}}\bm{y}^k_i,\;\;\; \bm{\til{y}}^k_i\leq u_{\bm{\alpha}_{k}}\bm{y}^k_i    &\quad \forall k \in \mcl{K}, \forall i \in \mcl{N}_{\bm{y}} &\\
 &\bm{\til{y}}^k_i \leq \ell_{\bm{\alpha}_k}\bm{y}^k_i + \bm{\alpha_{k}} - \ell_{\bm{\alpha}_k}, \;\;\; \bm{\til{y}}^k_i \geq  u_{\bm{\alpha}_k}\bm{y}^k_i + \bm{\alpha_{k}}- u_{\bm{\alpha}_k}  &\quad \forall k \in \mcl{K}, \forall i \in \mcl{N}_{\bm{y}}. &
\end{aligned}
\end{equation}
Inequalities (\ref{McCormik_explicit}) are tight when at least one of the two variables involved takes value equal to its upper or lower bound, otherwise they represent, in general, a relaxation. However, since variables $\yvk{k}$ are binary, the McCormick inequalities are tight at every feasible solution, for each valid choice of $\ell_{\bm{\alpha}_k}$ and $u_{\bm{\alpha}_k},\; k \in \Kset$. Since $\bm\alpha\in\mathbb{R}^{K}_{+}$ and $\bm e^\top \bm\alpha = 1$, as in \cite{hanasusanto2015k,vayanos2020robust}, $\ell_{\bm{\alpha}_k}$ and $u_{\bm{\alpha}_k},\; k \in \Kset$ can be set equal to 0 and 1, respectively. Nevertheless, the resulting linear relaxation may be extremely weak, slowing down the convergence of branch-and-cut algorithms. This is shown in the following example in which we illustrate how at fractional solutions the McCormick inequalities do not significantly constrain variables $\tyvk{k}$,\; $k \in \Kset$.
\example\label{example:num1} Consider a 2-adaptability instance of problem (\ref{mod_bi_problem}) with one recourse decision $\yvk{} \in\{0,1\}$ at a fractional solution $\yvk{1} = \yvk{2} = 0.5$ and $\alDV_1=\alDV_2 = 0.5$. Using the bounds $\alDV_{k} \in [\ell_{\alDV_{k}}, u_{\alDV_{k}}]=[0,1],\; k \in \Kset = \{1,2\}$ the McCormick inequalities (\ref{McCormik_explicit}) imply that $\tyvk{1}\in[0,\;\yvk{1}] = [0,\;0.5]$ and $\tyvk{2}\in[0,\yvk{2}] = [0,0.5]$. These bounds are somewhat trivial and provide a very loose approximation to constraints $\tyvk{1} =\alDV_1\yvk{1}= 0.025 $ and $\tyvk{2} =\alDV_2\yvk{2} = 0.025$.\endexample 
Example~\ref{example:num1} motivates us to examine how to reduce the bounds on $\alDV$ as these can significantly strengthen the McCormick inequalities and drastically impact the convergence rate of the branch-and-cut algorithm. To this end, we define set $\mathcal{A}_{B}$ which provides strengthened bounds on $\alDV$. 
\begin{equation*}\label{alpha_bounds_set}
\mcl{A}_B  := \left\{ \alDV \in \rrr^K_+ :\; \alDV_{k} \in \left[\ell_{\alDV_k},  u_{\alDV_k}\right]\; \forall k \in \Kset,\; \bm{e}^\top\alDV = 1 \right\}
\end{equation*}
where
\begin{equation*}\label{alpha_bounds}
\ell_{\alDV_k}=\left\{
\begin{array}{ll}
\frac{1}{K}, &\;\;\; \text{if $k = 1$}\\
0, &\;\;\; \text{otherwise}
\end{array}
\right. \quad \forall k \in \mcl{K},\qquad 
u_{\alDV_k}=  \frac{1}{k}  \quad\forall k \in \mcl{K}.
\end{equation*}
The following theorem ensures that the addition of constraints $\mcl{A}_B$ does not result in a conservative approximation. 
\begin{theorem}\label{th_alpha_bounds}
At least one optimal solution of problem (\ref{mod_bi_problem}) satisfies $\alDV \in \mcl{A}_B$. 
\end{theorem}

\addtocounter{example}{-1}
\example[continued] \label{example:num2} We now consider $\alDV \in \mcl{A}_S \cap \mcl{A}_B$, i.e.,~$\alDV_1 \in [0.5, 1]$ and $\alDV_2 \in [0, 0.5]$. We can immediately verify that inequalities~\eqref{McCormik_explicit} are tight, i.e.,~they ensure $\tyvk{1} =\alDV_1\yvk{1}= 0.025 $ and $\tyvk{2} =\alDV_2\yvk{2} = 0.025$.\endexample 
The McCormick inequalities can also be applied to the bilinear term $\gammaDV{k} \circ \bm \discV$ by introducing extra variables $\til{\gammaDV{k}}_i = 
\gammaDV{k}_i\discV_i$ for each $k \in \Kset$ and $i \in \mcl{N}_{\discV}$. However, 
since the $\gammaDV{k}$ always appears in problem~\eqref{mod_bi_problem} multiplied by variables $\discV$, one can linearize the bilinear terms by imposing the following 
constraint 
\begin{equation*}\label{lin_gamma_bigM}
    \gammaDV{k} \leq M\discV,\;\;\gammaDV{k} \geq -M\discV \quad \forall k \in \mcl{K},
\end{equation*}
where $M$ is a sufficiently large constant and replacing $\gammaDV{k}\circ \bm \discV$ with $\gammaDV{k}$, for each $k\; \in \Kset$ in 
constraints \eqref{mod_bil:cnstr:balance}-\eqref{mod_bil:cnstr:balance_k}. 

The linearized counterpart of problem~\eqref{mod_bi_problem} can be expressed as
\begin{equation}\label{mod_bi_problem_lin}
\begin{array}{>{\displaystyle}r@{\quad}>{\displaystyle}l@{\quad}>{\displaystyle}l}
 \min\quad&  \bm{b}^\top\left(\bm{\beta} + \sum_{k \in \mcl{K}}\bm{\beta}^{k}\right) & \\
\st \quad &  \bm{w} \in \mcl{W},\; \bm{\til{y}}^{k} \in \mcl{Y},\;\bm{\til{w}}^{k} \in \left[0, \frac{1}{k}\right]^{N_{\discV}}, \tyvk{k} \in \left[0, \frac{1}{k}\right]^{N_{\bm y}}, k \in \Kset&\\
&\bm{\alpha} \in \mcl{A}_S \cap \mcl{A}_B,\;\bm{\beta} \in \rrr^{L}_{+},\; \bm{\beta}^{k} \in \rrr^{L}_{+},\;  \bm{\gamma}^{k} \in \rrr^{N_{w}},\;k \in \Kset  &\\
& \bm{A}^\top\bm{\beta} = \sum_{k \in \mcl{K}}\bm{\gamma}^{k}  & \\
 &\left.\begin{array}{>{\displaystyle}l@{\quad}>{\displaystyle}l}
 \bm{A}^\top\bm{\beta}^{k} +\bm{\gamma}^{k} =  \bm{C}\bm{\til{w} }^{k}+\bm{P}{\bm{\til{y}}^{k}} & \\
 \gammaDV{k} \leq M\discV,\;\;\gammaDV{k} \geq -M\discV&\\
\bm{\til{y}}^k  \geq \frac{\mathbb{I}[k = 1]\bm{y}^k}{K},\;\;\bm{\til{y}}^k \geq  \frac{\bm{y}^k}{k} + \left(\bm{\alpha_{k}}- \frac{1}{k}\right)\bm{e} &\\
 \bm{\til{y}}^k \leq  \frac{\mathbb{I}[k = 1]\bm{y}^k}{K} + \left(\bm{\alpha_{k}} - \frac{\mathbb{I}[k = 1]}{K}\right)\bm{e},\;\;\bm{\til{y}}^k\leq \frac{\bm{y}^k}{k}&\\
 \bm{\til{w}}^k  \geq \frac{\mathbb{I}[k = 1]\bm{w}}{K},\;\;\bm{\til{w}}^k \geq  \frac{\bm{w}}{k} + \left(\alDV_{k}- \frac{1}{k}\right)\bm{e} &\\
 \bm{\til{w}}^k \leq  \frac{\mathbb{I}[k = 1]\bm{w}}{K} + \left(\alDV_k - \frac{\mathbb{I}[k = 1]}{K}\right)\bm{e},\;\;\bm{\til{w}}^k\leq \frac{\bm{w}}{k}&
 \end{array}    \right\}& \forall\; k \in \Kset.
 \end{array}
\end{equation}
For the remainder of the paper we will refer to  problem~\eqref{mod_bi_problem_lin} as the strengthened $K$-adaptability problem.
\subsection{Optimistic Inequalities}\label{sss:lb_cut}
In this section, we introduce a family of valid inequalities based on optimistic realizations of the uncertainty parameters. These so called optimistic inequalities leverage the following observation, which sheds light on the relationship between    $\alDV$ and $\left\{\yvk{k}\right\}_{k \in \Kset}$ in problem~\eqref{mod_bi_problem}. 
{\observation\label{observation}
 For fixed $(\discV,\{\yvk{k}\}_{k \in \Kset})$ we have
    \begin{align}
         &  \displaystyle\max_{\bm{\bar\xi}\in\Xi} \;\min_{k\in\mathcal{K}}\;\max_{\bm\xi\in\Xi(\bm w,\bm{\bar\xi})} \; \bm \xi^\top\bm C\bm w + \bm\xi^\top \bm P\bm y^k\nonumber\\[2ex]
         =&\left[
         \begin{array}{lll}
         \min & \tau\\
          \st  &\tau\in\mathbb{R},\,  \bm{\bar \xi}\in\Xi,\, \bm \xi^k \in \Xi(\bm w,\bm{\bar \xi})&\forall k\in\mathcal{K}\\
         & \tau \leq \bm \xi^\top\bm C\bm w + \bm\xi^\top \bm P\bm y^k & \forall k\in\mathcal{K}\\
         \end{array}\label{dual}
         \right] 
    \end{align}
where the equality follows from  \cite[Lemma 1]{vayanos2020robust}. Dualizing \eqref{dual} and combining with the minimization over $(\discV,\{\yvk{k}\}_{k \in \Kset})$ results in problem~\eqref{mod_bi_problem}. For each $k\in\mathcal{K}$, variable $\bm \alpha_k$ in problem~\eqref{mod_bi_problem} correspond to the dual multipliers associated with constraint $\tau \leq \bm \xi^\top\bm C\bm w + \bm\xi^\top \bm P\bm y^k$. At optimality $\tau = \bm \xi^\top\bm C\bm w + \bm\xi^\top \bm P\bm y^k$ for some $k\in\mathcal{K}$, with $\alDV_k>0$ indicating that the constraint is binding hence $\bm{y}^k$ influences the optimal value of problem~\eqref{mod_bi_problem}, while $\alDV_k=0$ indicates the contrary. This is also evident in problem~\eqref{mod_bi_problem} where if $\alDV_k = 0$ for some $k\in\mathcal{K}$, then $\bm{y}^k$ does not influence the optimal value since term $\bm{\alpha}_k\bm{y}^k$ in constraint~\eqref{mod_bil:cnstr:balance_k} is responsible for the interaction of $\bm{y}^k$ with the rest of the variables. Hence for $k\in\mathcal{K}$ with $\bm\alpha_k = 0$ setting $\bm y^k := \bm y^{\prime k}$ to any $\bm y^{\prime  k}\in\mathcal{Y}$, then  $(\bm w, \{\bm y_k\}_{k\in\mathcal{K}})$ achieves the same optimal value in \eqref{mod_bi_problem}.
}

We now propose the following inequalities which are based on an optimistic view of uncertainty. Let $\bm{\zeta}^{\discV} \in \rrr^{N_{\discV}},\   \bm{\zeta}^{\yvk{}} \in \rrr^{N_{\bm{y}}}$ be two vectors such that $\bm{\zeta}^{\bm{w}}_{i} \leq {\bm{\xi}}^{\top}\bm{C}_i$ for  $i \in \mcl{N}_{\discV}$  and  $\bm{\zeta}^{\bm y}_{j} \leq {\bm{\xi}}^{\top}\bm{P}_j $ for  $j \in \mcl{N}_{\bm{y}}$ for all  $\uncV \in \uncSet$, where $\bm{C}_i$ and $\bm{P}_j$  represent the $i$-th and $j$-th columns of matrices $\bm{C}$ and $\bm{P}$, respectively. Then the following family of inequalities
\begin{equation}
\bm{b}^\top\left(\bm{\beta} + \sum_{k \in \mcl{K}}\bm{\beta}^{k}\right) 
\geq {\bm{\zeta}^{\discV}}^\top\discV+ {\bm{\zeta}^{\yvk{}}}^\top\yvk{k}  \quad \forall k \in \mcl{K}, \label{best_scenario_cuts}
\end{equation}
can be added to problem~\eqref{mod_bi_problem_lin} without loss of optimality, as stated in the following theorem.
\begin{theorem}
At least one optimal solution of problem~\eqref{mod_bi_problem_lin} satisfies inequalities~\eqref{best_scenario_cuts}. \label{th_optimistic_cut}
\end{theorem}
\noindent
In general, there may be many possible vectors $\bm{\zeta}^{\discV},\; \bm{\zeta}^{\yvk{}}$, however, the tightest vectors for which inequalities are still valid can be generated by solving $N_{\discV}+N_{\bm{y}}$ linear optimization problems
\begin{equation*}
\bm{\zeta}^{\discV}_{i} = \min_{\uncV \in \uncSet}  \uncV^\top \bm{C}_{i}\quad \forall i \in \mcl{N}^{\discV},\qquad\bm{\zeta}^{\yvk{}}_{j} = \min_{\uncV \in \uncSet} \uncV^\top \bm{P}_{j} \quad \forall j \in \mcl{N}^{\yvk{}}
\end{equation*}
The way vectors $\bm{\zeta}^{\discV}$ and $\bm{\zeta}^{\yvk{}}$ are generated suggests that the effectiveness of (\ref{best_scenario_cuts}) strongly depends on the considered uncertainty set. Note that a similar inequality can be added to problem~\eqref{original_prob}. In this case, the inequality only involves $\discV$ and take the form $\displaystyle \Phi(\discV) \geq {\bm{\zeta}^{\discV}}^\top\discV$.

\subsection{Reformulation-Linearization-Technique Inequalities}\label{sss:rlt_cuts}
The \emph{Reformulation-Linearization-Technique} (RLT), introduced in \cite{sherali1992new}, is a procedure used to derive tight relaxations for discrete and continuous non-convex problems. It consists of two  phases: reformulation and linearization. In the reformulation phase, valid quadratic constraints are generated using pairwise product operations between constraints and/or variables. In the linearization phase, the products resulting from the reformulation phase are linearized through additional variables and McCormick inequalities.

We apply the RLT scheme to problem~\eqref{mod_bi_problem_lin} to obtain tighter linear relaxations that can assist the convergence of the branch-and-cut algorithm. In our implementation, the reformulation phase consists in generating valid inequalities by multiplying each constraints in the definition of $\ySet$ with $\alDV_k$ for each $k \in \Kset$. Notice that if some constraints in the definition of $\ySet$ are added on the fly, additional RLT inequalities can be generated within the branch-and-cut procedure. This is the case for the orienteering problem whose constraints are described in \eqref{Y_set}. In this case  the exponential number of subtours elimination constraints are dynamically added when needed together with the corresponding RLT inequalities.

\section{Computational Experiments}\label{section:comp}
In this section, we present experiments to assess the performance of the exact algorithm presented in Section~\ref{sec:exa_algo} using the information inequalities \eqref{exact_information_cut}, and the strengthened $K$-adaptability problem~\eqref{mod_bi_problem_lin} in conjunction with optimistic inequalities and the RLT inequalities presented in Section \ref{sec:k_adapt}. Section~\ref{section:öp_experiments} compares the performance for the two solution schemes on the sensor placement orienteering problem. As the strengthened $K$-adaptability formulation can also be applied to problem instances without decision dependent information discovery, in  Section~\ref{section:shortest_path} we benchmark its performance to the $K$-adaptability formulation proposed by \cite{hanasusanto2015k} and the iterative solution scheme proposed by \cite{subramanyam2019k} using the shortest path problem, as both publications present results on this problem instance. All of our experiments are performed in single core on a machine equipped with a 4.0GHz Intel i7-600K processor and with 24GB of RAM, using CPLEX 12.10. A time of 7200 seconds is allowed to solve each instance.
\subsection{Computational Experiments for the Sensor Placement Orienteering Problem}\label{section:öp_experiments}
In the following we use the orienteering problem networks  presented in \cite{tsiligirides1984heuristic} which consists of three network structures which we denote by TS1, TS2 and TS3, involving 30, 19 and 31 nodes, respectively. For all networks, the start and destination nodes are placed in the center of the of the network. We also construct smaller versions of the same networks where we remove half of the nodes closest to the start/destination node, resulting in networks with 15, 10 and 16 nodes, respectively. The original networks are denoted by TS1N30, TS2N19 and TS3N31, and their smaller variants by TS1N15, TS2N10 and TS1N16, respectively. Similar to \cite{tsiligirides1984heuristic}, due to the different topology of each network, different values of the maximum time duration $T$ are considered, and are summarized in Table~\ref{tab:op_instances}. All instances can be found in~\cite{instances}.

We consider the sensor placement orienteering problem~\eqref{SensorPlacement} where $\mathcal{Y}$ is given in \eqref{Y_set} and the uncertainty set is given by $\Xi = \{\bm \xi\in[0,U]^N : \bm e^\top \bm\xi = 1\}$ which expresses the view that the total profit in the network is known and equals to 1 (or 100\%), however, the geographical breakdown of the profit in each node is uncertain. The upper bound $U$ denotes the maximum profit of each node, and its value for each network is given in Table \ref{tab:op_instances}. To assess the impact of the number of sensors on the total profit collect, we denote by $\delta$ the fraction of the total nodes where sensors can be placed, and for each $\delta \in\{0.25,0.5,0.75\}$ the admissible sensor placements are described by $\displaystyle\discVdom = \left\{\discV \in \bbb^{N_{\discV}} : \bm{e}^\top\discV \leq \lceil \delta N_{\discV}\rceil\right\}$.
We next compare the relative performance of the exact method and the strengthened $K$-adaptability formulation (Section~\ref{sec:exa_vs_kada}), and examine the behavior of  the strengthened $K$-adaptability formulation (Section~\ref{section:performance_kada}).

\subsubsection{Comparing the Exact Algorithm and \texorpdfstring{$K$}{}-adaptability }\label{sec:exa_vs_kada}
In our first experiment, we compare the performance of the exact solution method and of the strengthened $K$-adaptability formulation for $K\in\{2,3,4\}$ in terms of their optimality gap. For each $\delta = \{0.25,0.5,0.75\}$ and network topology, we solve using both methods problem \eqref{SensorPlacement} for all duration times $T$ listed in Table \ref{tab:op_instances}.  Notice that the optimal value of problem \eqref{exact_master}, which we denote by ``best progressive bound", provides a progressive bound to the optimal value of problem \eqref{SensorPlacement}. Hence, the optimality gap after 7200 seconds of computational time, is calculated as  $\frac{\text{``best conservative bound" - ``best progressive bound"}}{\text{``best progressive  bound}"}$ where the ``best conservative bound" is the best objective value achieved by each method within the time limit. Table~\ref{tab:exa_vs_k_ada} reports the average gap with respect to the different values of $T$ achieved for each method. Tables~\ref{tab:exa} and \ref{tab:results:perf:kada} (next section)  presents further details about the behavior of exact method and strengthened $K$-adaptability, respectively. 

\begin{table}[htbp!]
\begin{footnotesize}
\begin{center}
  \caption{Optimality gap for the exact solution method and the $K$-adaptability approximation.}
          \sffamily
  \renewcommand{\arraystretch}{0.8}
  
  \begin{tabular}{cc}
      \begin{tabular}[t]{|lc|ccc|c|}
    \hline
    \hline
    Network & $\delta$ & 2-Adapt & 3-Adapt & 4-Adapt & Exact \\
   \hline
        \multirow{3}[1]{*}{TS1N15} & 0.25  & 8.9\%   & 7.1\%    & 7.1\%    & 0.0\%  \\
          & 0.50   & 16.2\%   & 8.8\%    & 7.4\%    & 0.0\%  \\
          & 0.75  & 17.7\%   & 10.2\%   & 7.6\%    & 0.0\%  \\
          \hline
    \multirow{3}[1]{*}{TS2N10} & 0.25  & 11.1\%   & 0.0\%    & 0.0\%    & 0.0\%  \\
          & 0.50   & 16.8\%   & 1.8\%    & 0.0\%    & 0.0\%  \\
          & 0.75  & 17.0\%   & 4.4\%    & 0.2\%    & 0.0\%  \\
 \hline
    \multirow{3}[2]{*}{TS3N16} & 0.25  & 7.1\%    & 2.4\%    & 0.0\%    & 0.0\%  \\
          & 0.50   & 12.9\%   & 5.4\%    & 2.9\%    & 0.0\%  \\
          & 0.75  & 13.9\%   & 6.4\%    & 3.9\%    & 0.0\%  \\
    \Xhline{1pt}
    \end{tabular}%
    & \hfill  
    \begin{tabular}[t]{|lc|ccc|c|}
    \hline
    \hline
    Network & $\delta$ & 2-Adapt & 3-Adapt & 4-Adapt & Exact \\
   \hline
    \multirow{3}[2]{*}{TS1N30} & 0.25 & 10.0\%  & 8.0\%   & 11.2\%  & 9.5\%  \\
          & 0.50 & 8.6\%   & 4.5\%   & 4.8\%   & 3.2\%  \\
          & 0.75 & 8.6\%   & 3.7\%   & 2.8\%   & 0.3\%  \\
    \hline
  \multirow{3}[2]{*}{TS2N19} & 0.25 & 8.6\%   & 7.6\%   & 7.6\%   & 7.5\%  \\
          & 0.50 & 21.0\%  & 8.3\%   & 2.4\%   & 0.6\%  \\
          & 0.75 & 22.0\%  & 10.5\%  & 5.1\%   & 0.0\%  \\
     \hline
   \multirow{3}[2]{*}{TS3N31} & 0.25 & 9.1\%   & 7.8\%   & 11.2\%  & 10.8\%  \\
          & 0.50 & 8.2\%   & 4.4\%   & 3.5\%   & 3.3\%  \\
          & 0.75 & 8.2\%   & 3.8\%   & 2.0\%   & 0.8\%  \\
    \Xhline{1pt}
    \end{tabular}%
    \\
  \end{tabular}
  \label{tab:exa_vs_k_ada}%
  \end{center}
  \end{footnotesize}
\end{table}%

\begin{table}[htbp!]
\begin{footnotesize}
\begin{center}
  \caption{Additional details on the performance of the exact algorithm. We denote by \emph{Opt (\#)} the number of instances solved to optimality, by \emph{Time (s)} the average computational time (in seconds) for instances solved to optimality and by \emph{Gap} the average optimality gap for the instances not solved within the time limit. Overall, 182 instances out of 258 are solved to optimality with an average computational time of 411 seconds. The average gap for the 76 instances not solved to optimality is equal to 7.94\%.}
          \sffamily
  \renewcommand{\arraystretch}{0.8}
  
  \begin{tabular}{cc}
      \begin{tabular}[t]{|lc|ccc|}
    \hline
    \hline
    Network & $\delta$ & Opt (\#) & Time (s) & Gap  \\
   \hline
        \multirow{3}[1]{*}{TS1N15} & 0.25  & 14/14   & 80.4   & -    \\
         & 0.50  & 14/14  &  89.5  &  -  \\
         & 0.75  & 14/14  & 12.4  &  -  \\
          \hline
    \multirow{3}[1]{*}{TS2N10} & 0.25  &9/9    & 3.8  &    -  \\
          & 0.50   & 9/9    &  4.1  &-    \\
          & 0.75  & 9/9 &  1.6  &    - \\
 \hline
    \multirow{3}[2]{*}{TS3N16} & 0.25  &14/14 &     43.1&  -   \\
          & 0.50   & 14/14   &   156.7  & -  \\
          & 0.75  &  14/14  & 15.4   &    - \\
    \Xhline{1pt}
    \end{tabular}%
    & \hfill  
    \begin{tabular}[t]{|lc|ccc|}
    \hline
    \hline
 Network & $\delta$ & Opt (\#) & Time (s) & Gap  \\
   \hline
    \multirow{3}[2]{*}{TS1N30} & 0.25 & 6/18  & 71.2  & 14.0\%   \\
          & 0.50 & 6/18     &  103.0 &   4.8\%  \\
          & 0.75 &  10/18   &   193.2 &  0.7\%  \\
    \hline
  \multirow{3}[2]{*}{TS2N19} & 0.25 & 7/11   & 2364.4   & 17.1\%   \\
          & 0.50 & 11/11   & 2985.1  &  -   \\
          & 0.75 & 11/11  & 1155.9 &   - \\
     \hline
   \multirow{3}[2]{*}{TS3N31} & 0.25 & 6/20  & 154.6   & 14.5\%   \\
          & 0.50 &  7/20 & 321.2 &   6.6\% \\
          & 0.75 &  7/20  & 120.6  &     1.2\%\\
    \Xhline{1pt}
    \end{tabular}%
    \\
  \end{tabular}
  \label{tab:exa}%
  \end{center}
  \end{footnotesize}
\end{table}%

From our first experiment we observed the following. First notice that the exact method optimally solves the smaller instances TS1N15, TS2N10 and TS2N16 for all values of $\delta$ and $T$. Overall, it solves to optimality 70\% of the total instances. 
Interestingly, for the larger more challenging instances, we observe that the exact method is more efficient in finding near optimal solutions as $\delta$ increases.  The strengthened $K$-adaptabililty approximation is solved to optimality for the majority of the smaller instances as well, while across all values of  $K\in\{2,3,4\}$ around 65\% are solved to optimality, see Table \ref{tab:results:perf:kada}.  Hence, for the smaller instances, the gap endured from the $K$-adaptability approximation results primarily from the lack of adaptivity from the recourse decisions $\yvk{}$, which gets worst as the  $\delta$ increases. Still, as expected, the gap improves as the level of adaptability $K$ increases and in some cases $K$-adaptability achieves the optimal solution of problem \eqref{SensorPlacement}, e.g., for network TS2N10 and $\delta = 0.25$, $K = 3$ achieves the optimal solution for all values of $T$. We conclude the following: $(i)$ the exact method provides on average better quality solutions compared to $K$-adaptability, however, in two cases (i.e., TS1N30 and TS3N31 with $\delta = 0.25$) $K$-adaptability can provide better solutions  with small values of $K$ regardless if the optimization problems are not solved to optimality, (ii) if adaptability is limited, the extra information provided by increasing the number of sensors, i.e., increasing $\delta$, will not be able to be exploited by $K$-adaptability. We will see further evidence of this in the following and in Section~\ref{sec:case}.

Note that our results do not report 1-adaptability as this case is not interesting. This is the case as there is essentially no adaptability, thus increasing $\delta$ will not improve the optimal solution. In fact all sensor placements will achieve the same optimal value as the case where no sensors are used. Similar to Example~\ref{example:sol:two_stage_prob}, due to the structure of uncertainty set $\Xi$, the worst-case $\uncV$ will move the profit away from the visited nodes. Hence, for small values of $T$, 1-adaptability will achieve zero profit regardless of which nodes are visited, while for larger values of $T$ it will achieve the minimum profit allowable by the uncertainty set. In terms of computational time, 1-adaptability can be optimally solved within a few seconds.

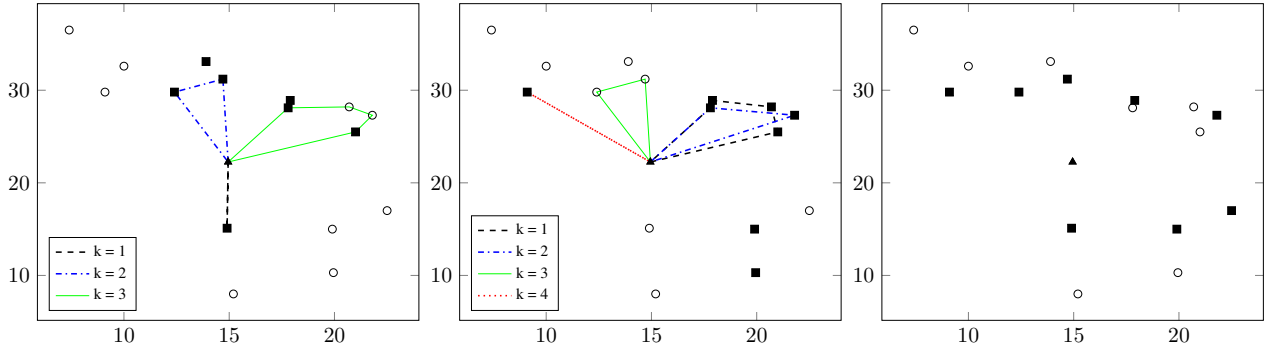
\begin{figure}[!htbp]
    \captionsetup[subfigure]{font=footnotesize,labelfont=footnotesize}
\begin{subfigure}[t]{0.32\textwidth}
\centering
\begin{tikzpicture}[scale=0.74]
\begin{axis} [
    legend pos=south west,
scatter/classes={%
	start={mark=triangle*},
	end={mark=triangle},
discovery_on={mark=square*},
discovery_off={mark=o}}]
\addplot[scatter, only marks,  scatter src = explicit symbolic, forget plot] table[meta = label]{
x     y      label
14.950000	22.250000	start
14.950000	22.250000	end
20.700000	28.200000	discovery_off
21.800000	27.300000	discovery_off
22.500000	17.000000	discovery_off
19.900000	15.000000	discovery_off
14.900000	15.100000	discovery_on
12.400000	29.800000	discovery_on
17.800000	28.100000	discovery_on
9.100000	29.800000	discovery_off
10.000000	32.600000	discovery_off
13.900000	33.100000	discovery_on
19.950000	10.300000	discovery_off
15.200000	8.000000	discovery_off
14.700000	31.200000	discovery_on
7.400000	36.500000	discovery_off
21.000000	25.500000	discovery_on
17.900000	28.900000	discovery_on
};
	\addplot[thick, dashed, mark=none] coordinates {
(14.95, 22.25)
(14.9, 15.1)
(14.95, 22.25)
};
	\addplot[thick,dash dot,mark=none, color=blue] coordinates {
(14.95, 22.25)
(14.7, 31.2)
(12.4, 29.8)
(14.95, 22.25)
};
	\addplot[mark=none, color=green] coordinates {
(14.95, 22.25)
(17.8, 28.1)
(20.7, 28.2)
(21.8, 27.3)
(21, 25.5)
(14.95, 22.25)
};
\addlegendentry{{\scriptsize k = 1}}
\addlegendentry{{\scriptsize k = 2}}
\addlegendentry{{\scriptsize k = 3}}
\addlegendentry{{\scriptsize k = 4}}
\end{axis}
\end{tikzpicture}
\caption{3-adaptability: worst-case profit 3\%.}\label{fig:op_sol:result:3adapt}
\end{subfigure}
\hfill
\begin{subfigure}[t]{0.32\textwidth}
\centering
\begin{tikzpicture}[scale=0.74]
\begin{axis} [legend pos=south west,
scatter/classes={%
	start={mark=triangle*},
	end={mark=triangle},
discovery_on={mark=square*},
discovery_off={mark=o}}]
\addplot[scatter, only marks,  scatter src = explicit symbolic, forget plot] table[meta = label]{
x     y      label
14.950000	22.250000	start
14.950000	22.250000	end
20.700000	28.200000	discovery_on
21.800000	27.300000	discovery_on
22.500000	17.000000	discovery_off
19.900000	15.000000	discovery_on
14.900000	15.100000	discovery_off
12.400000	29.800000	discovery_off
17.800000	28.100000	discovery_on
9.100000	29.800000	discovery_on
10.000000	32.600000	discovery_off
13.900000	33.100000	discovery_off
19.950000	10.300000	discovery_on
15.200000	8.000000	discovery_off
14.700000	31.200000	discovery_off
7.400000	36.500000	discovery_off
21.000000	25.500000	discovery_on
17.900000	28.900000	discovery_on
};
	\addplot[thick, dashed, mark=none] coordinates {
(14.95, 22.25)
(21, 25.5)
(20.7, 28.2)
(17.9, 28.9)
(17.8, 28.1)
(14.95, 22.25)
};
	\addplot[thick,dash dot,mark=none, color=blue] coordinates {
(14.95, 22.25)
(17.8, 28.1)
(21.8, 27.3)
(14.95, 22.25)
};
	\addplot[mark=none, color=green] coordinates {
(14.95, 22.25)
(12.4, 29.8)
(14.7, 31.2)
(14.95, 22.25)
};
	\addplot[thick,dotted,color=red,mark=none] coordinates {
(14.95, 22.25)
(9.1, 29.8)
(14.95, 22.25)
};
\addlegendentry{{\scriptsize k = 1}}
\addlegendentry{{\scriptsize k = 2}}
\addlegendentry{{\scriptsize k = 3}}
\addlegendentry{{\scriptsize k = 4}}
\end{axis}
\end{tikzpicture}
\caption{4-adaptability: worst-case profit 5\%.}\label{fig:op_sol:result:4adapt}
\end{subfigure}
\hfill
\begin{subfigure}[t]{0.32\textwidth}
\begin{tikzpicture}[scale= 0.74]
\begin{axis} [legend pos=south west,
scatter/classes={%
	start={mark=triangle*},
	end={mark=triangle},
discovery_on={mark=square*},
discovery_off={mark=o}}]
\addplot[scatter, only marks,  scatter src = explicit symbolic, forget plot] table[meta = label]{
x     y      label
14.950000	22.250000	start
14.950000	22.250000	end
20.700000	28.200000	discovery_off
21.800000	27.300000	discovery_on
22.500000	17.000000	discovery_on
19.900000	15.000000	discovery_on
14.900000	15.100000	discovery_on
12.400000	29.800000	discovery_on
17.800000	28.100000	discovery_off
9.100000	29.800000	discovery_on
10.000000	32.600000	discovery_off
13.900000	33.100000	discovery_off
19.950000	10.300000	discovery_off
15.200000	8.000000	discovery_off
14.700000	31.200000	discovery_on
7.400000	36.500000	discovery_off
21.000000	25.500000	discovery_off
17.900000	28.900000	discovery_on
};
\end{axis}
\end{tikzpicture}
\caption{Optimal: worst-case profit 7.7\%.}
\label{fig:op_sol:result:exa}
\end{subfigure}
\caption{Sensor placements for TS3N16 with $\delta = 0.50$ and $ T = 20$ generated by $K$-adaptability ($K = 3, 4$) and the exact method. Nodes are plotted in cartesian coordinates using black squares (sensor placed) or  circles (sensor not placed).
}
\label{fig:op_sol:result}
\end{figure}

Next we discuss the sensor placement positions resulting from $K$-adaptability for $K \in\{3,4\}$ and the exact solution method. We consider the network TS3N16 with $\delta = 0.5$ and time duration $T=20$, i.e., the network has 16 nodes and at most 8 sensors can be placed throughout the network. All problems are solved to optimality. Figure~\ref{fig:op_sol:result} presents the nodes of the network (white circles) and the optimal sensor positions for each method (black squares). The start/destination node is denoted by a black triangle. For $K$-adaptability the graph also presents the resulting $K$ constant decisions $\yvk{}$. Note that we do not report results for 2-adaptability  as the optimal profit is zero. This implies that regardless of the sensor placement, due to not enough adaptability from the recourse decisions, nature is able to adversary place the worst-case profit in other nodes. In contrast, 3- and 4-adaptability produce positive worst-case profit which amounts to 3\% and 5\% of the total profit in the network, see  Figures~\ref{fig:op_sol:result:3adapt}~and~\ref{fig:op_sol:result:4adapt}, respectively. Similar to 2-adaptability, 3-adaptability is not able to fully utilize the information from all 8 sensors, as the optimal solution places only 7 out of the 8 sensors. In contrast, 4-adaptability places all 8 sensors and achieves a higher collected profit. The sensor placement for both 3- and 4-adaptability is significantly different to the sensor placement of the exact solution, Figure \ref{fig:op_sol:result:exa}, which places the sensors almost uniformly around the start/destination node, and   is able to achieve a worst-case profit of 7.7\%. 

To better understand the difference in exact and $K$-adaptability approximation, and to distinguish if the suboptimality is due to the suboptimal placement of the sensors or the lack of adaptivity from the $K$-adaptable solution, we next evaluate $\Phi(\discV^{*K})$ using Algorithm \ref{Algo_CCG},  where $\bm \discV^{*3}$ and $\bm \discV^{*4}$ are the optimal sensor placement  resulting from 3- and 4-adaptability, respectively.  We observe that $\Phi(\bm \discV^{*3})$ results in a 5\% profit, and $\Phi(\bm \discV^{*4})$ results in in 6.67\% profit, as opposed to the 3\% and 5\% indicated by the optimal value of problem~ \eqref{mod_bi_problem_lin}, respectively. We conclude that indeed, part of the suboptimality $K$-adaptability is due to the lack of adaptability, and part is due suboptimal placement of the sensors. 

\subsubsection{Further Details on the Performance of the Strengthened \texorpdfstring{$K$}{}-adaptability Approximation}\label{section:performance_kada}
We next take a closed look at the behavior of the strengthened $K$-adaptability formulation and compare its performance to the 
$K$-adaptability formulation proposed in \cite{vayanos2020robust}. Results 
are reported in Table~\ref{tab:results:perf:kada}. 
\emph{Str $K$-adapt} denotes the strengthened formulation introduced in Section~\ref{sec:k_adapt} and  \emph{$K$-adapt} refer to the formulation proposed in~\cite{vayanos2020robust}.
The last column labeled \emph{RN} reports the average root node relaxation improvement of the strengthened formulation with respect to original one (calculated as $\frac{\text{Root Node strengthened $K$-adapt -  Root Node $K$-adapt}}{\text {Root Node strengthened $K$-adapt}}$). The rest of the table is interpreted as Table~\ref{tab:exa}. We clarify that in Table~\ref{tab:results:perf:kada}, \emph{Gap} refers to the optimality gap of the mixed-integer solver when the time limit is reached. All the average values are calculated with regards to the different values of the  maximum time durations $T$.  Overall, the strengthened formulation consistently outperforms the formulation of \cite{vayanos2020robust}, solving to optimality 65\% of all instances as opposed to 52\% of the \cite{vayanos2020robust} formulation. Beside that, optimality gaps are consistently reduced when instances are not solved to optimality. The benefit of the improved bounds and valid inequalities  can be also noticed by observing the bound at the root node relaxation of the two formulation. As shown in the last column, the improvements of the root node relaxation ranges between 45\% and 80\%. 

In the majority of cases, the stronger lower bounds assist the branch-and-cut procedure to converge faster. However, there is a small number of exceptions, e.g.,~one instance of TS1N15 with $\delta = 0.25$ and $K = 2$, where the strengthened formulation has an adverse effects. Nevertheless, we can conclude that for the sensor placement orienteering problem, the strengthened $K$-adaptability formulation improves significantly the computational time.
\begin{table}[!htbp]
\begin{scriptsize}
\begin{center}
  \caption{Summary of the results for the $K$-adaptability formulations.}
\centering
        \sffamily
 \renewcommand{\arraystretch}{1.15}
    \begin{tabular}{|lcc|cc|cc|cc|c|}
    \hline
    \hline
      &       &       & \multicolumn{2}{c|}{Opt (\#)} & \multicolumn{2}{c|}{Time (s)} & \multicolumn{2}{c|}{Gap} &  \multicolumn{1}{c|}{} \\
    Network & $K$     & $\delta$ & Str $K$-Adapt & $K$-Adapt & Str $K$-Adapt      & $K$-Adapt      & Str $K$-Adapt & $K$-Adapt & RN  \\
    \hline
    \multirow{9}[5]{*}{TS1N15} & \multirow{3}[2]{*}{2} & 0.25  & 13/14    & 14/14    & 324.3 & 114.8 & 7.1\%   & -   & 58.5\%  \\
          &       & 0.50   & 14/14    & 14/14    & 15.4  & 50.5  & -   & -   & 56.3\%   \\
          &       & 0.75  & 14/14    & 14/14    & 14.4  & 44.4  & -   & -  & 53.7\%  \\
\cline{2-10}          & \multirow{3}[2]{*}{3} & 0.25  & 9/14     & 8/14     & 289.9 & 167.9 & 19.2\%  & 24.1\%  & 64.5\%  \\
          &       & 0.50   & 14/14    & 10/14    & 213.4 & 426.9 & -   & 11.2\%  & 64.5\%  \\
          &       & 0.75  & 14/14    & 10/14    & 94.1  & 505.7 & -   & 11.0\%  & 64.8\%  \\
\cline{2-10}          & \multirow{3}[1]{*}{4} & 0.25  & 5/14     & 6/14     & 246.0 & 462.4 & 21.4\%  & 38.0\%  & 76.1\%   \\
          &       & 0.50   & 9/14     & 7/14     & 552.4 & 316.0 & 4.8\%   & 27.6\%  & 76.2\%  \\
          &       & 0.75  & 11/14    & 6/14     & 1792.3 & 114.3 & 4.3\%   & 23.6\%  & 75.1\%  \\
          \hline
    \multirow{9}[5]{*}{TS2N10} & \multirow{3}[1]{*}{2} & 0.25  & 9/9     & 9/9     & 0.9   & 0.3   & -   & -   & 42.6\%  \\
          &       & 0.50   & 9/9     & 9/9     & 0.5   & 0.2   & -  & -   & 43.1\%  \\
          &       & 0.75  & 9/9     & 9/9     & 0.6   & 0.2   & -   & -   & 42.3\%   \\
\cline{2-10}          & \multirow{3}[2]{*}{3} & 0.25  & 9/9     & 9/9     & 30.3  & 6.1   & -   & -   & 53.5\% \\
          &       & 0.50   & 9/9     & 9/9     & 3.6   & 3.8   & -   & -   & 52.9\%  \\
          &       & 0.75  & 9/9     & 9/9     & 4.6   & 3.7   & -   & -   & 53.4\%   \\
\cline{2-10}          & \multirow{3}[2]{*}{4} & 0.25  & 7/9     & 9/9     & 357.8 & 789.8 & 7.7\%   & -   & 59.2\%   \\
          &       & 0.50   & 9/9     & 9/9     & 33.7  & 131.4 & -   & -   & 63.3\%  \\
          &       & 0.75  & 9/9     & 9/9     & 43.2  & 73.2  & -   & -   & 62.4\%   \\
    \hline
    \multirow{9}[6]{*}{TS3N16} & \multirow{3}[2]{*}{2} & 0.25  & 14/14    & 14/14    & 21.5  & 8.8   & -   & -   & 45.9\%   \\
          &       & 0.50   & 14/14    & 14/14    & 6.3   & 9.5   & -   & -   & 46.9\%  \\
          &       & 0.75  & 14/14    & 14/14    & 6.2   & 9.7   & -   & -   & 46.3\%   \\
\cline{2-10}          & \multirow{3}[2]{*}{3} & 0.25  & 13/14    & 11/14    & 392.2 & 715.9 & 12.0\%  & 12.4\%  & 61.1\%  \\
          &       & 0.50   & 14/14    & 12/14    & 274.8 & 1164.3 & -   & 10.2\%  & 60.9\%  \\
          &       & 0.75  & 14/14    & 13/14    & 197.5 & 592.4 & -   & 13.8\%  & 60.9\%   \\
\cline{2-10}          & \multirow{3}[2]{*}{4} & 0.25  & 7/14     & 4/14     & 312.5 & 402.4 & 16.1\%  & 22.4\%  & 67.1\%   \\
          &       & 0.50   & 10/14    & 6/14     & 1068.2 & 1460.0 & 7.1\%   & 20.3\%  & 66.5\%  \\
          &       & 0.75  & 11/14    & 6/14     & 968.6 & 552.6 & 6.8\%   & 16.7\%  & 68.3\%  \\
    \hline
    \multirow{9}[5]{*}{TS1N30} & \multirow{3}[2]{*}{2} & 0.25  & 12/18    & 6/18     & 675.1 & 95.6  & 22.2\%  & 21.5\%  & 49.4\%  \\
          &       & 0.50   & 15/18    & 6/18     & 415.7 & 101.7 & 2.0\%   & 17.3\%  & 49.2\%   \\
          &       & 0.75  & 16/18    & 8/18     & 404.3 & 516.7 & 1.7\%   & 17.8\%  & 49.4\%  \\
\cline{2-10}          & \multirow{3}[2]{*}{3} & 0.25  & 6/18     & 4/18    & 1232.2 & 168.1 & 13.5\%  & 36.3\%  & 71.3\%   \\
          &       & 0.50   & 5/18     & 4/18     & 1275.7 & 455.9 & 5.9\%   & 36.5\% & 70.2\%   \\
          &       & 0.75  & 6/18     & 4/18     & 665.9 & 411.3 & 4.6\%  & 35.3\% & 71.2\%   \\
\cline{2-10}          & \multirow{3}[1]{*}{4} & 0.25  & 4/18     & 3/18     & 481.3 & 0.1   & 20.3\% & 43.2\%  & 77.5\%  \\
          &       & 0.50   & 3/18     & 3/18    & 0.1   & 0.5   & 11.5\%  & 41.4\%  & 80.4\%   \\
          &       & 0.75  & 3/18     & 3/18     & 0.2   & 0.3   & 10.0\% & 42.0\%  & 79.9\%   \\
          \hline
    \multirow{9}[4]{*}{TS2N19} & \multirow{3}[1]{*}{2} & 0.25  & 9/11     & 11/11    & 219.6 & 391.7 & 3.7\%  &-  &49.6\%    \\
          &       & 0.50   & 11/11    & 11/11   & 162.1 & 179.2 & -  & -  & 50.1\%  \\
          &       & 0.75  & 11/11    & 11/11    & 151.1 & 324.8 & -  & -   & 49.4\%  \\
\cline{2-10}          & \multicolumn{1}{c}{\multirow{3}[2]{*}{3}} & 0.25  & 1/11     & 3/11     & 4222.3 & 1024.4 & 21.3\% & 35.3\%  & 60.8\%   \\
          &       & 0.50   & 8/11     & 6/11     & 649.4 & 666.1 & 2.7\%  & 13.4\% & 60.8\%   \\
          &       & 0.75  & 8/11     & 7/11     & 998.5 & 792.2 & 1.4\%  & 16.1\%  & 60.8\%  \\
\cline{2-10}          & \multirow{3}[1]{*}{4} & 0.25  & 1/11     & 1/11     & 0.3   & 0.0   & 32.9\% & 42.4\%  & 65.9\%  \\
          &       & 0.50   & 2/11     & 2/11     & 1895.1 & 273.8 & 5.9\%   & 28.0\%  & 65.5\%  \\
          &       & 0.75  & 3/11     & 2/11     & 3115.4 & 614.6 & 2.1\%   & 24.1\%  & 66.7\%   \\
          \hline
    \multirow{9}[5]{*}{TS3N31} & \multicolumn{1}{c}{\multirow{3}[1]{*}{2}} & 0.25  & 15/20    & 7/20    & 547.4 & 38.4  & 2.9\%  & 21.8\%  & 48.0\%   \\
          &       & 0.50   & 18/20    & 8/20     & 338.6 & 209.5 & 4.4\%  & 18.9\% & 48.1\%   \\
          &       & 0.75  & 18/20    & 8/20     & 487.0 & 60.8  & 2.5\%   & 19.2\% & 48.5\%   \\
\cline{2-10}          & \multicolumn{1}{c}{\multirow{3}[2]{*}{3}} & 0.25  & 4/20    & 4/20     & 2000.7 & 34.5  & 10.5\%  & 34.3\%  & 67.3\%  \\
          &       & 0.50   & 10/20    & 4/20     & 2048.8 & 79.8  & 7.2 \%  & 32.6\% & 67.1\%  \\
          &       & 0.75  & 6/20     & 4/20     & 1513.3 & 37.8  & 4.0 \%  & 32.5\%  & 68.6\%  \\
\cline{2-10}          & \multirow{3}[2]{*}{4} & 0.25  & 5/20     & 3/20     & 1197.6 & 0.2   & 20.8\%  & 39.9\%  & 74.2\% \\
          &       & 0.50   & 4/20     & 4/20     & 328.3 & 438.1 & 9.9\%   & 38.1\%  & 74.0\%  \\
          &       & 0.75  & 5/20     & 4/20     & 927.0 & 393.2 & 8.4\%   & 37.3\%  & 73.1\%  \\
    \hline
    Total   &       &       & 502/774 & 405/774 &       &       &       &       &        \\
    Average   &       &       &       &       & 359.4 & 275.4 & 11.5\%  & 30.0\%  & 61.8\%   \\
    \Xhline{1pt}
  \end{tabular}%
   \label{tab:results:perf:kada}%
\end{center}
\end{scriptsize}\end{table}%
\subsection{Assessing the Strengthened \texorpdfstring{$K$}{}-adaptability Formulation: The Shortest Path Problem}\label{section:shortest_path}
The strengthened $K$-adaptability formulation is also valid for problem instances without decision dependent information discovery. To assess its relative merits compared to other approaches in the literature, we consider the robust shortest path problem solved in \cite{hanasusanto2015k} and in \cite{subramanyam2019k}. Note that this problem is not characterized by decision dependent information discovery: uncertainty is always revealed before the selection of the recourse decisions. In our formulation, this is achieved by fixing $\discV = \bm{e}$.
As in \cite{hanasusanto2015k} and in \cite{subramanyam2019k} we consider a shortest path problem defined on a directed weighted graph $\mcl{G} = (\mcl{N}, \mcl{A})$, where $\mcl{N} := \{1, \ldots N\}$ represents the set of the $N$ nodes, $\mcl{A}  \subseteq \mcl{N} \times \mcl{N}$ is the set of arcs and $\bm{c}_{ij}(\uncV) =(1+\uncV_{ij}/2)\til{\bm{c}}_{ij}$ is the cost associated to arc $(i,j)$, function of the Euclidean distance $\til{\bm{c}}_{ij}$ and of the uncertain vector $\uncV$ which belongs to the uncertainty set $\uncSet \subseteq \rrr^{|\mcl{A}|}$. The objective is to determine $K$ paths starting from a start node $s \in \mcl{N}$ and terminating at a terminal node $t \in \mcl{N} $, $t \neq s$ before observing $\uncV$, such that the worst-case length of the best path among the $K$ defined is minimized. The problem can be formulated as follows.
\begin{equation*}\label{shortest_path_formulation}
    \begin{array}{r@{\quad}l}
\displaystyle\min &\displaystyle \max_{\uncV \in \uncSet}\; \min_{k \in \Kset} \sum_{(i,j) \in A}\bm{c}_{ij}(\uncV)\yvk{k}_{ij}  \\
\displaystyle \text{\st}&\displaystyle \yvk{k}_{ij} \in \bbb,\; \forall (i,j) \in A, \forall k \in \Kset \\
        &\displaystyle \sum_{(j,l) \in A}\yvk{k}_{jl} \geq \sum_{(i,j) \in A}\yvk{k}_{ij}+\mathbb{I}[j = s] - \mathbb{I}[j = t]  \quad \forall j \in V,\ \forall k \in \Kset
    \end{array}
\end{equation*}
 We consider instances with $N\in\{30, 40, 50\}$ using the same generation procedure as in  \cite{hanasusanto2015k} and in \cite{subramanyam2019k}. In each problem, the location of each node is chosen uniformly at random from the interval $[0, 10]^2$. Start and terminal nodes are chosen as the pair of nodes with largest Euclidean distance. The set of arc $\mcl{A}$ is generated by removing from the set $\mcl{N} \times \mcl{N}$ the $\left \lfloor0.7|\mcl{N} \times \mcl{N}|\right \rfloor$ arcs with the largest nominal weight. A budget uncertainty set  $\uncSet = \{\uncV \in [0,1]^{A} : \sum\limits_{(i,j) \in A} \uncV_{ij} \leq B   \}$ set is considered. 
For each graph size, 10 problems are generated and solved for $B$ equal to 3 and 6. 

 Results are reported in Table~\ref{tab:shortest_path_comparison} and the interpretation of its columns is the same as in Tables~\ref{tab:exa} and \ref{tab:results:perf:kada}. 
\begin{table}[htbp]
\begin{footnotesize}
\begin{center}
    \caption{Computational results for the shortest path problem}
    \sffamily
\renewcommand{\arraystretch}{0.8}
    \begin{tabular}{|lcc|ccc|ccc|ccc|}
    \hline
    \hline
     &   &       & \multicolumn{3}{c|}{\cite{hanasusanto2015k}} & \multicolumn{3}{c|}{\cite{subramanyam2019k}} & \multicolumn{3}{c|}{Str $K$-Adapt} \\
  &  $K$     & $N$     & Opt (\#)   & Time (s)  & Gap  & Opt (\#)   & Time (s)  & Gap & Opt (\#)   & Time (s)  & Gap  \\
 \hline
 \multirow{9}[2]{*}{$B$ = 3}  &   \multirow{3}[2]{*}{2} & 30    & 0/10     & -     & 49.2\%  & 7/10     & 852.0   & 6.7\%   & 10/10    & 2.7   & - \\
 &         & 40    & 0/10     & -     & 63.7\%  & 1/10     & 1034.0  & 7.0\%   & 10/10    & 20.0  & - \\
 &         & 50    & 0/10     & -     & 69.5\%  & 0/10     & -     & 9.3\%   & 10/10    & 158.1 & - \\
     \cline{2-12}
 &   \multirow{3}[2]{*}{3} & 30    & 0/10     & -     & 58.8\%  & 2/10     & 987.2 & 3.8\%   & 10/10    & 8.4   & - \\
 &         & 40    & 0/10     & -     & 70.0\%  & 0/10     & -     & 6.5\%   & 10/10    & 60.3  & - \\
 &         & 50    & 0/10     & -     & 74.8\%  & 0/10     & -     & 8.4\%   & 10/10    & 1207.0 & - \\
     \cline{2-12}
 &   \multirow{3}[2]{*}{4} & 30    & 0/10     & -     & 70.4\%  & 2/10     & 1875.7 & 4.7\%   & 10/10    & 51.0  & - \\
  &        & 40    & 0/10     & -     & 76.1\%  & 0/10     & -     & 6.1\%   & 10/10    & 232.7 & - \\
  &        & 50    & 0/10     & -     & 82.2\%  & 0/10     & -     & 7.3\%   & 10/10    & 2451.4 & - \\
    \hline
 \multirow{9}[2]{*}{$B$ = 6}  & \multirow{3}[2]{*}{2} & 30    & 0/10     &   -    & 50.4\%  & 0/10     &  -     & 10.8\%  & 10/10    & 27.5\%  & - \\
  &        & 40    & 0/10     &   -    & 67.1\%  & 0/10     &    -   & 15.4\%  & 10/10    & 321.6 & - \\
 &         & 50    & 0/10     &  - & 68.4\%  & 0/10     &   -    & 17.2\%  & 9/10     & 1915.2 & 0.4\%  \\
   \cline{2-12}
&    \multirow{3}[2]{*}{3} & 30    & 0/10     &     -  & 60.0\%  & 0/10     &   -    & 10.8\%  & 10/10    & 113.9 & - \\
 &         & 40    & 0/10     &  -     & 72.7\%  & 0/10     &      - & 14.1\%  & 9/10     & 939.2 & 0.6\% \\
 &         & 50    & 0/10     &   -    & 77.0\%  & 0/10     &  -     & 15.4\%  & 4/10     & 2714.0 & 3.4\% \\
    \cline{2-12}
 &   \multirow{3}[2]{*}{4} & 30    & 0/10     &  -     & 70.7\%  & 0/10     &     -  & 9.5\%   & 10/10    & 529.3 & - \\
  &        & 40    & 0/10     &    -   & 79.2\%  & 0/10     &      - & 12.3\%  & 6/10     & 2324.9 & 1.3\% \\
  &        & 50    & 0/10     &     -  & 86.4\% & 0/10     &      - & 13.3\%  & 3/10     & 2528.5 & 3.6\% \\
    \Xhline{1pt}
    \end{tabular}%
  \label{tab:shortest_path_comparison}%
  \end{center}
  \end{footnotesize}
\end{table}%
The results clearly demonstrate the merits of the strengthened  formulation which provides significant computational benefits compared to the existing methods.

\section{Case Study: Collecting Medicine Crates at the Alrijne Hospital}\label{sec:case}
The Alrijne hospital in the Netherlands is a general hospital with more than 3700 staff and has a capacity of 500 beds across several departments, treating yearly more than 190,000 patients. In this study we consider the hospital in Leiderdorp which is a single building with six floors. The medications are supplied by a pharmacy operating from a depot at the basement of the hospital.
To ensure the quality of the medicines, the deliveries are done in specialized crates equipped with sensors that record temperature, light and shock variations to make sure that medication is not damaged and that the crates are not opened by unauthorized personnel. The crates are expensive equipment that need to be recovered to be reused in future deliveries. 

The pharmacy is facing the problem that although the crates are delivered without problems to their destinations, their collection process is highly problematic as a large number of these crates are lost. Consequently, to achieve timely delivery of the medicine, the pharmacy is forced to acquire more of these specialized crates, thus increasing  significantly their operating costs.  In an attempt to understand this issue, in collaboration with the Alrijne hospital, we conducted a pilot study from February 11 to April 30 in 2019 (\cite{alrijneMasterThesis,alrijneBachelorThesis}). The primary purpose of the study was to understand where the crates ended up. To this end, transmitters were retrofitted in the crates and their location was tracked by a network of sensors placed in a part of the hospital.  The study disputed the original belief that the crates where stolen and moved out of the hospital, and in fact showed that the crates were ending up in different departments within the hospital.  

In this follow up study, we integrate the placing of sensors with the crate collection process. Following  the guidelines set by the hospital, the aim is to be able to efficiently collect the crates by placing a limited number of sensors throughout the hospital, thus limiting the infrastructural as well as the amount of sensory data that need to be collected, stored and analyzed. The aim is to understand the trade-off between the number of sensors, the flexibility of adapting the collection routes and the maximum available duration of the collection routes. 
\subsection{Data Description}
\begin{figure}[ht]
 \captionsetup[subfigure]{font=footnotesize,labelfont=footnotesize}
\begin{subfigure}[b]{.45\linewidth}
\centering
\includegraphics[width=0.8\textwidth]{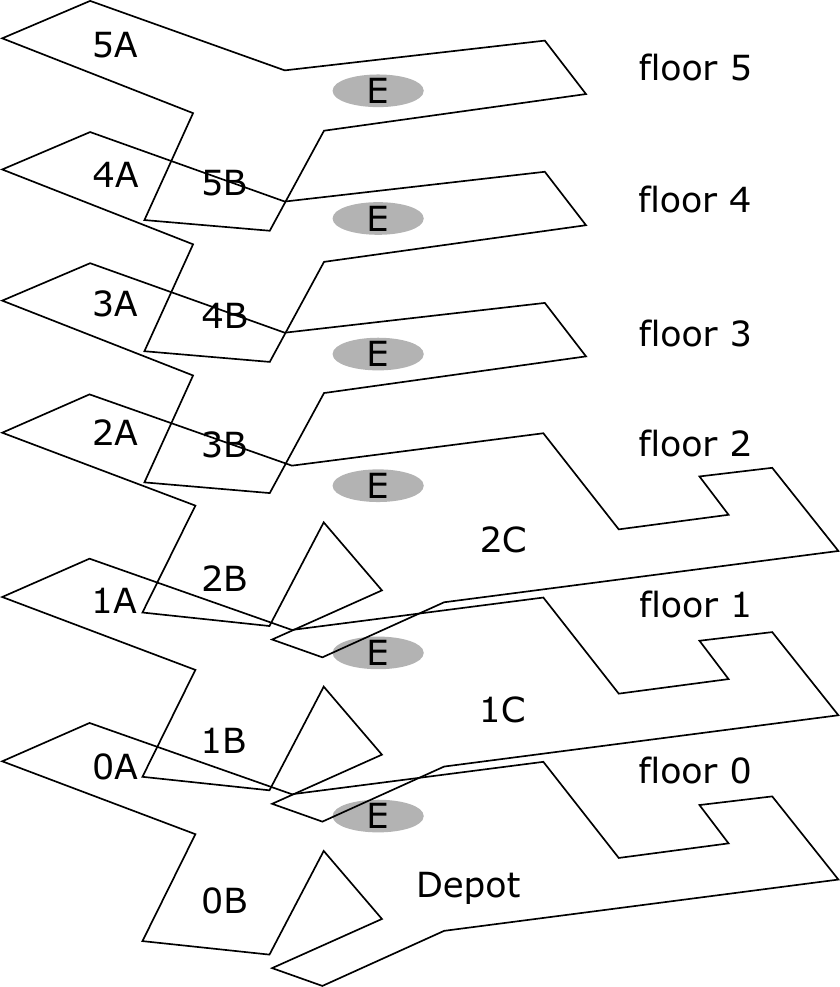}
\caption{Alrijne hospital floors map.}\label{subfig:Alrijne:floormap}
\end{subfigure}
\begin{subfigure}[b]{.45\linewidth}
\centering
\begin{tikzpicture}[->,>=stealth',shorten >=2pt, node distance=1.5cm,
                    semithick,framed,background rectangle/.style={thick,fill=gray1,draw=black}]]
    \tikzstyle{every state}=[fill=white,draw=black,text=black,font=\small]

\node[ellipse,draw] (Elv) at (0,0) {E};

 \node[triangle, draw] (dep) at (-1.5,2) {\scriptsize{D}};
 \node[shape=circle,draw,above of = dep] (n0B) {0B};
 \node[shape=circle,draw, left of = dep] (n0C) {0C};

\draw[-,shorten >= -.5pt]  (dep) to node [auto] {} (n0C); 
\draw[-,shorten >= -.5pt,shorten <= -.5pt]  (n0C) to node [auto] {} (n0B); 
\draw[-,shorten >= -.5pt,shorten <= -.5pt]  (dep) to node [auto] {} (n0B); 
\draw[-,shorten >= -.5pt,shorten <= -.5pt]  (Elv) to node [auto] {} (n0B); 
\draw[-,shorten >= -.5pt,shorten <= -.5pt]  (Elv) to node [auto] {} (n0C); 
\draw[-,shorten >= -.5pt,shorten <= -.5pt]  (Elv) to node [auto] {} (dep);

\node[shape=circle,draw] (n1A) at (1.5,2.5) {1A};
\node[shape=circle,right of=n1A,draw] (n1C) {1C};
 \node[shape=circle, draw, above left of = n1A] (n1B) {1B};
 
\draw[-,shorten >= -.5pt,shorten <= -.5pt](n1B)  to node [auto] {} (n1A); 
\draw[-,shorten >= -.5pt,shorten <= -.5pt,bend left](n1B)  to node [auto] {} (n1C); 
\draw[-,shorten >= -.5pt,shorten <= -.5pt](n1C)  to node [auto] {} (n1A); 
\draw[-,shorten >= -.5pt,shorten <= -.5pt]  (Elv) to node [auto] {} (n1A); 
\draw[-,shorten >= -.5pt,shorten <= -.5pt]  (Elv) to node [auto] {} (n1C); 
\draw[-,shorten >= -.5pt,shorten <= -.5pt]  (Elv) to node [auto] {} (n1B); 

\node[shape=circle,draw] (n2C) at (2,0) {2C};
 \node[shape=circle,draw, below right of = n2C] (n2A) {2A};
 \node[shape=circle, draw, above right of = n2C] (n2B) {2B};

\draw[-,shorten >= -.5pt,shorten <= -.5pt](n2B)  to node [auto] {} (n2A); 
\draw[-,shorten >= -.5pt,shorten <= -.5pt](n2B)  to node [auto] {} (n2C); 
\draw[-,shorten >= -.5pt,shorten <= -.5pt](n2A)  to node [auto] {} (n2C); 
\draw[-,shorten >= -.5pt,shorten <= -.5pt]  (Elv) to node [auto] {} (n2A); 
\draw[-,shorten >= -.5pt,shorten <= -.5pt]  (Elv) to node [auto] {} (n2C); 
\draw[-,shorten >= -.5pt,shorten <= -.5pt]  (Elv) to node [auto] {} (n2B);

\node[shape=circle,draw] (n3A) at (1,-2.5) {3A};
 \node[shape=circle,draw, right of = n3A] (n3B) {3B};
 
\draw[-,shorten >= -.5pt,shorten <= -.5pt] (n3B)  to node [auto] {} (n3A); 
\draw[-,shorten >= -.5pt,shorten <= -.5pt]   (Elv) to node [auto] {} (n3A); 
\draw[-,shorten >= -.5pt,shorten <= -.5pt]   (Elv) to node [auto] {} (n3B);

\node[shape=circle,draw] (n4A) at (-2,-2.5) {4A};
 \node[shape=circle,draw, right of = n4A] (n4B) {4B};
 
\draw[-,shorten >= -.5pt,shorten <= -.5pt] (n4B)  to node [auto] {} (n4A); 
\draw[-,shorten >= -.5pt,shorten <= -.5pt]  (Elv) to node [auto] {} (n4A); 
\draw[-,shorten >= -.5pt,shorten <= -.5pt]   (Elv) to node [auto] {} (n4B);

\node[shape=circle,draw] (n5A) at (-3,0.5) {5A};
\node[shape=circle,draw, below of=n5A] (n5B) {5B};
\draw[-,shorten >= -.5pt,shorten <= -.5pt] (n5B)  to node [auto] {} (n5A); 
\draw[-,shorten >= -.5pt,shorten <= -.5pt]   (Elv) to node [auto] {} (n5A); 
\draw[-,shorten >= -.5pt,shorten <= -.5pt]   (Elv) to node [auto] {} (n5B); 

\end{tikzpicture} 

\caption{Graph representation of Alrijne hospital.}\label{subfig:Alrijne:graph}
\end{subfigure}
\caption{The Alrijne hospital floor map and the graph representing connections between the different areas and floors of the hospital. The 27 departments are clustered into 14 areas (Figure~\ref{subfig:Alrijne:floormap}). The departments in each cluster are physically very close to each other. As illustrated in Figure \ref{subfig:Alrijne:floormap}, each floor is divided in at most three areas with area A and B corresponding to the north and south wing of the hospital respectively and area C covers the right wing and the south zone. Node D (triangle) represents the pharmacy depot and node E the elevators. Arcs between nodes indicates the direct path between two areas.}\label{fig:Alrijne}
\end{figure}
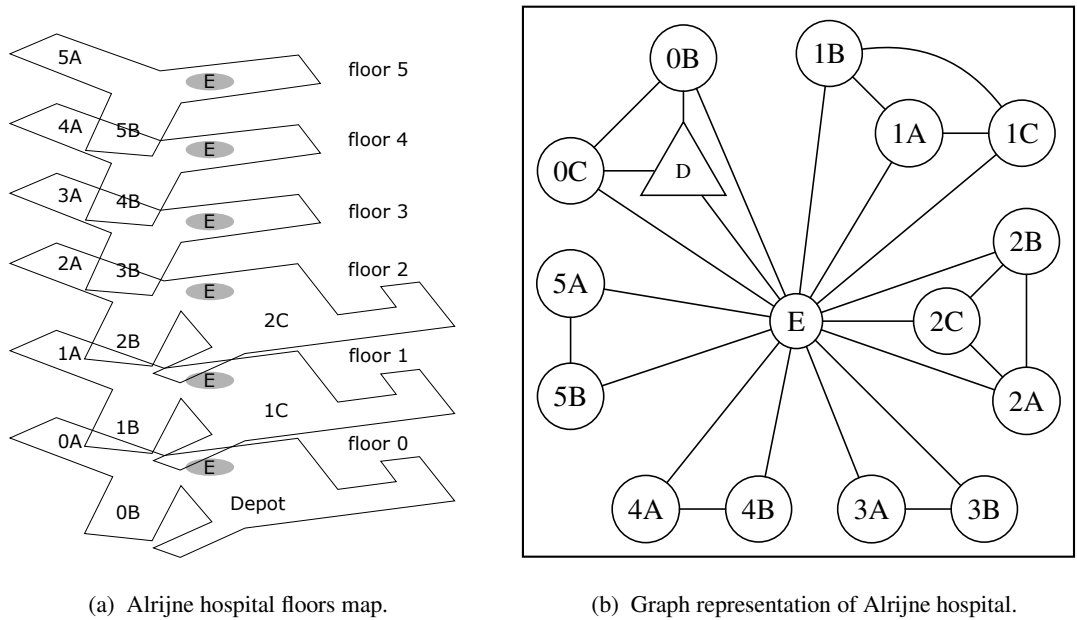

We model the problem as a sensor placement orienteering problem~\eqref{SensorPlacement}. The hospital is represented via a graph where the 27 departments of the hospital are clustered into 14 nodes that cover the 6 floors of the building. The clustering was based on the floor plan of the hospital and ensures that all departments within each node have short walking distances, see Table~\ref{tab:Alrijne:departments} and Figure~\ref{fig:Alrijne}. The associated travel time $t_{(i,j)}$ between nodes $i$ and $j$ is calculated as the estimated shortest path walking time plus the time need to collect the crates in node $j$. If two nodes are located in different floors of the hospital, the travel time is calculated as the walking time from the departing node to the elevator, the time spent in the elevator and the walking time from the elevator to the destination node plus the time to collect the crates. The resulting distance matrix is available at~\cite{instances}. We model by $\bm \xi_i$ as the uncertain percentage of crates located in node $i$. As before, the uncertainty set is given by $\Xi_{1} = \{\bm \xi \in [0,U]^{14} :  \bm e^\top \bm{\xi} = 1\}$ with $U = 0.2$ and reflects that the total number of crates circulating in the hospital is known through constraints $\bm e^\top \bm{\xi} = 1$,
however, their geographical breakdown is uncertainty, with at most 20\% of the total crates can be found at any node of the network.

\subsection{Adaptability and Value of Information}
\begin{figure}[!ht]
    \centering
    \resizebox{0.4\textwidth}{!}{%
\begin{tikzpicture}
\begin{axis}[
    xlabel={$\delta$},
    ylabel={Worst-case profit (\%)},
    xmin=0, xmax=1.20,
    ymin=55, ymax=90,
    xtick={.0,.25,.50,.75,1.00},
    ytick={60,70,80,90,100},
    legend pos=north west,
    ymajorgrids=true,
    grid style=dashed,
]
\addplot[
    color=black,
    mark=*,
    ]
    coordinates {
    (0,60)(.25,66.7)(.50,77.78)(.75,84.6)(1.00,85.7)
    };
    \label{alrijne_plot_exa}
    \addplot[
    color=black,
    mark=square,
    ]
    coordinates {
    (.0,60)(.25,60)(.50,60)(.75,60)(1.00,60)
    };
    \label{alrijne_plot_2ad}
     \addplot[
    color=black,
    mark=triangle,
    ]
    coordinates {
    (.0,60)(.25,66.7)(.50,66.7)(.75,66.7)(1.00,66.7)
    };
    \label{alrijne_plot_3ad}
       \addplot[
    color=black,
    mark=o,
    ]
    coordinates {
    (0,60)(.25,66.7)(.50,75.0)(.75,75.0)(1.00,75.0)
    };
    \label{alrijne_plot_4ad}
\end{axis}
\node [draw,fill=white] at (rel axis cs: 0.15,0.8) {\shortstack[l]{
\ref{alrijne_plot_exa} {\scriptsize Exact} \\
\ref{alrijne_plot_2ad} {\scriptsize 2-adapt} \\
\ref{alrijne_plot_3ad} {\scriptsize 3-adapt }\\
\ref{alrijne_plot_4ad} {\scriptsize 4-adapt}
}};
\end{tikzpicture}
}
    \caption{Worst-case profit of the 2-, 3-, 4-adaptable and exact  solution method for increasing values of $\delta$ and for $T=\text{2 hours}$.}
    \label{fig:Alrijne:exp:gain}
\end{figure}
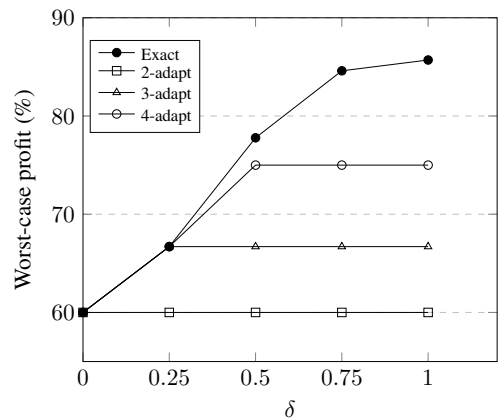

In our first experiment, we assess the impact of the number of sensors placed in the network on the percentage of crates collected, using the exact algorithm and $K$-adaptability with $K \in \{2,3,4\}$. As before, we denote use  $\discVdom = \{\discV \in \bbb^{N_{\discV}} : \bm{e}^\top\discV \leq \lceil \delta N_w\rceil\}$ and examine the impact of the number of sensors by changing $\delta\in \{0,\ 0.25,\ 0.50,\ 0.75,\ 1\}$. We set the maximum time duration to two hours, i.e., $T = 2$~hours. All problem instances were solved to optimality. Figure~\ref{fig:Alrijne:exp:gain} reports the optimal value achieved by both methods as a function of $\delta$. As expected, the exact algorithm achieves higher worst-case profit, and is able to better utilize the information from the increased number of sensors. We again observe, that 
$K$-adaptability can only utilize a proportion of the information provided. Indeed, 2-adaptability approximation seems extremely restrictive and achieves the same solution as $\delta = 0$, while 3-adaptability and 4-adaptability can effectively utilize the information from 25\% and 50\% of the sensors, respectively. Despite this shortcoming, a salient feature of $K$-adaptability is that the optimal solution provides $K$ alternative routing plans from which one can be selected once measurements from the sensors has been obtained. This could have a practical interest if the personnel needs to train for standardized routing. However, due to the relatively small size of the hospital, once the $\ov{\bm{\xi}}$ associated to node equipped with sensor has been observed, the optimal worst-case routing problem can be achieved by solving $\displaystyle \max_{\yvk{} \in \ySet} \min_{\bm{\xi}\in\Xi(\discV,\ov{\bm{\xi}})} \uncV^{\top}\yvk{}$. From the experiment we conclude that extra information can only be assimilated if the recourse decisions have adequate adaptability.

\subsubsection{Impact of the Maximum Route Duration and Uncertainty Set}
In our second experiment, we examine the impact of the maximum time duration $T$ on the quality of the solution. To this end, we again evaluate the solution of the exact algorithm for different the values of $\delta$ by changing the time duration $T \in \{0.5, 1, 1.5, 2\}$ hours and considering the same uncertainty set adopted in the first experiment. The results are presented in Figure~\ref{subfig:Alrijne:comparison_T_US_u}. As expected, we observe that as the time duration increases, the collection process improves as the routing plan is able to visit more nodes in the network. As a result the impact of uncertainty diminishes. In contrast, the information provided by sensors is more impactful for short time durations. We conclude that as the hospital is trying to limit its investment in sensor infrastructure, they should consider the trade-off between number of sensors and the duration of the collecting process.
 \begin{figure}[!ht]
  \captionsetup[subfigure]{font=footnotesize,labelfont=footnotesize}
\begin{subfigure}[b]{.48\linewidth}
\centering
\resizebox{!}{0.70\textwidth}{%
\begin{tikzpicture}
\begin{axis}[
    xlabel={$\delta$},
    ylabel={Worst-case profit (\%)},
    xmin=0, xmax=1.20,
    ymin=0, ymax=100,
    xtick={0,.25,.50,.75,1.00},
    ytick={0,20,40,60,80,100},
    legend pos=north west,
    ymajorgrids=true,
    grid style=dashed,
]
\addplot[
    color=gray,
    mark=triangle,
    ]
    coordinates {
    (0,0)(.25,0)(.50,0)(.75,12)(1.00,14)
    };
    \label{alrijne_plot_1800}
    \addplot[
    color=gray,
    mark=square,
    ]
    coordinates {
    (0,0)(.25,10)(.50,33)(.75,38)(1.00,38)
    };
    \label{alrijne_plot_3600}
     \addplot[
    color=black,
    mark=triangle*,
    ]
    coordinates {
    (0,0)(.25,40)(.50,58)(.75,64)(1.00,64)
    };
    \label{alrijne_plot_5400}
       \addplot[
    color=black,
    mark=square*,
    ]
    coordinates {
    (0,60)(.25,67)(.50,78)(.75,85)(1.00,86)
    };
    \label{alrijne_plot_7200}
\end{axis}
\node [draw,fill=white] at (rel axis cs: 0.17,0.88) {\shortstack[l]{
\ref{alrijne_plot_7200} {\scriptsize $T = 2\;\text{hrs}$}\\
\ref{alrijne_plot_5400} {\scriptsize $T = 1.5\;\text{hrs}$}\\
\ref{alrijne_plot_7200} {\scriptsize $T = 1\;\text{hrs}$}\\
\ref{alrijne_plot_1800} {\scriptsize $T = 0.5\;\text{hrs}$}
}};
\end{tikzpicture}
}
\caption{Worst-case profit using $\Xi_1$.}\label{subfig:Alrijne:comparison_T_US_u}
\end{subfigure}
\hfill
\begin{subfigure}[b]{.48\linewidth}
\centering
\resizebox{!}{0.70\textwidth}{%
\begin{tikzpicture}
\begin{axis}[
    xlabel={$\delta$},
    xmin=0, xmax=1.20,
    ymin=0, ymax=100,
    xtick={0,.25,.50,.75,1.00},
    ytick={20,40,60,80,100},
    legend pos=north west,
    ymajorgrids=true,
    grid style=dashed,
]
\addplot[
    color=gray,
    mark=triangle,
    ]
    coordinates {
    (0,5)(.25,7)(.50,12)(.75,14)(1.00,14)
    };
    \addplot[
    color=gray,
    mark=square,
    ]
    coordinates {
    (0,12)(.25,34)(.50,38)(.75,39)(1.00,39)
    };
     \addplot[
    color=black,
    mark=triangle*,
    ]
    coordinates {
    (0,46)(.25,52)(.50,60)(.75,64)(1.00,64)
    };
       \addplot[
    color=black,
    mark=square*,
    ]
    coordinates {
    (0,82)(.25,82)(.50,82)(.75,85)(1.00,86)
    };
\end{axis}
\end{tikzpicture}
}
\caption{Worst-case profit using $\Xi_2$.}\label{subfig:Alrijne:comparison_T_US_ui}
\end{subfigure}
\caption{\mbox{Comparison between  optimal worst-case profits obtained for different values of $\delta$ and $T$ using different uncertainty sets.}}\label{fig:Alrijne:comparison_T_US}
\end{figure}
  
In our last experiment, we examine the situation where additional information is known on the number of crates at each node. 
In particular, we use data collected during the pilot  to estimate the expected fraction of crates $\widehat{U}_i$ at each node $i$ based on the orders of each department. 
We consider the uncertainty set $\Xi_{2} = \left\{\bm{\xi}\in\rrr^{14}_{+}: \bm{\xi}_i\in\left[\widehat{U}_i\left(1-\theta\right) ,\widehat{U}_i\left(1 +\theta\right)\right],\, i \in \{1, \ldots 14\},\; \bm{e}^\top \bm{\xi} = 1\right\}$ with $\theta = 0.75$.
Table \ref{tab:Alrijne:departments} summarizes the values of $\widehat{U}_i$ for each node.  We solve the problem using the exact algorithm for different values of $\delta$ and $T$. The results are presented in Figure \ref{subfig:Alrijne:comparison_T_US_ui}. We observe that the sensor information does provide additional benefits for all values of $T$, however, their impact is less significant compared to the case where uncertainty set $\Xi_1$ is considered. This can be explained by the fact that $\Xi_2$ ensures that a minimum amount of crates is present at each node  and that for some nodes the maximum value of $\widehat{U}_i$ is quite low, e.g., $\widehat{U}_3 = 0.05$, therefore due to the presence of constraint $\bm e^\top \bm\xi = 1$, the orienteering routing plan will focus on nodes which have higher potential of collecting crates and avoid visiting nodes which have smaller values of $\widehat{U}_i$. There are two take away from this experiment. First, if additional information is available which limits the amount of uncertainty in the model, then the impact of the number of sensors diminishes compared to the case where uncertainty is uniform across all nodes, making the nodes somewhat ``indistinguishable". Second, from a practical point of view, since problem (\ref{SensorPlacement}) tries to model effectively an infinite horizon problem where if crates are not collected at a given node then they could potentially accumulate over time, having a constant $U$ across all nodes ensures that the sensors are not biasedly placed in the network.

\ACKNOWLEDGMENT{
The authors gratefully acknowledge the financial support by the Dutch Institute of Advanced Logistics (Dinalog) [grant number 2017-2-131TKI].}

\bibliography{bibliography}
\bibliographystyle{abbrvnat}



\newpage

\ECSwitch



~\vspace{-1.5cm}

\ECHead{Electronic Companion}

\section{Proofs for Section~\ref{sec:exa_algo}}

\textsc{Proof of Proposition~\ref{preposition_tight_and_valid}.} To prove the proposition, we need to show that the two following properties hold.
\begin{enumerate}
\item[$(i)$] $\Phi(\discV) = \Phi(\discV) - \left(\Phi(\discV) - \un{\Phi}\right)\phi(\discV, \discV)$,\label{int_cut_prop_1}
\item[$(ii)$] $\Phi(\discV) \geq \Phi(\discV') - \left(\Phi(\discV') - \un{\Phi}\right)\phi(\discV', \discV)\quad \forall \discV' \in \discVdom$.\label{int_cut_prop_2}
\end{enumerate}
Property $(i)$ holds since $\phi(\discV, \discV) = 0$. This property ensures tightness (i.e., the right hand side of the inequality coincides with the value of $\Phi(\discV)$ at point $\discV$). Moreover, for $\discV' \in \discVdom \setminus \{\discV\}$, we observe that $\phi(\discV', \discV)\geq 1$. As $\un{\Phi} \leq \Phi(\discV)$ for all $\discV \in \discVdom$, property $(ii)$  ensures that the inequalities never over estimate the value of $\Phi(\discV)$, thus the inequalities are valid.
\Halmos\endproof\vspace{5pt}

\noindent \textsc{Proof of Theorem~\ref{th_exa_converge}.} 
The finite converge of Algorithm~\ref{Algo_Logic_Benders} relies on $\discVdom$ being a finite set and on the convexity of the lower envelope of $\Phi(\discV)$. Since by Proposition~\ref{preposition_tight_and_valid} the inequalities~\eqref{exact_convex_envelope} are valid and tight, in the worst-case the algorithm will generate all possible $\discV\in \discVdom$ thus terminating in $|\discVdom|$ iterations.
\Halmos\endproof\vspace{5pt}

 \noindent\textsc{Proof of Lemma~\ref{lemma_exa_function_prop}.}
Assume that there exist at least one $i \in \mcl{N}_w$ such that $\discV^1_i > \discV^2_i$.
By definition, $\uncSet(\discV^1, \ov{\uncV})\ \subseteq\ \uncSet(\discV^2, \ov{\uncV})$ and therefore $\displaystyle \max_{\bm{\ov\xi} \in \uncSet} \min_{\bm y\in\mathcal{Y}} \max_{\bm \xi\in\Xi(\bm{w}^{1}, \bm{\ov \xi})}\bm \xi^\top \bm{P}\bm{y} \leq \max_{\bm{\bar\xi} \in \uncSet} \min_{\bm y\in\mathcal{Y}} \max_{\bm \xi\in\Xi(\discV^{2}, \bm{\ov\xi})}\bm \xi^\top \bm{P}\bm{y}$ which implies $\Phi(\discV^1) \leq \Phi(\discV^2)$, as desired.
\Halmos\endproof\vspace{5pt}

\noindent\textsc{Proof of Proposition~\ref{proposition_information_inequalities}.}
To prove the proposition, we need to show the following conditions hold:
\begin{enumerate}
\item[$(i)$] $\Phi(\discV) = \Phi(\discV) - \left(\Phi(\discV) - \un{\Phi}\right)\rho(\discV, \discV)$,\label{int_lift_cut_prop_1}
\item[$(ii)$] $\Phi(\discV) \geq \Phi(\discV') - \left(\Phi(\discV') - \un{\Phi}\right)\rho(\discV', \discV)\quad \forall\; \discV' \in \discVdom$,\label{int_lift_cut_prop_2}
\item[$(iii)$] $\Phi(\discV') - \left(\Phi(\discV') - \un{\Phi}\right)\rho(\discV', \discV'')\geq \Phi(\discV') - \left(\Phi(\discV') - \un{\Phi}\right)\phi(\discV', \discV'' )\quad \forall\; \discV'\; \in\; \discVdom$.\label{int_lift_cut_prop_3}
\end{enumerate}
Condition $(i)$ ensures tightness at integer solutions and holds since $\rho(\discV, \discV) = 0$.
To show that condition $(ii)$ holds we need to distinguish two cases. If $\discV > \discV'$, the cut is valid since  $\rho(\discV', \discV) \geq 1$. When $\discV \leq \discV'$, $\rho(\discV', \discV) =0$ 
that implies $\Phi(\discV) \geq \Phi(\discV')$ which follows from Lemma \ref{lemma_exa_function_prop}.
Condition $(iii)$ ensures that the information inequalities are stronger than inequalities~\eqref{exact_convex_envelope} at fractional solutions. As $\Phi(\discV) \geq \un{\Phi}$, it is sufficient to prove that $\rho(\discV, \discV') \leq \phi(\discV, \discV') $. 
We have $\phi(\discV, \discV') - \rho(\discV, \discV') = \sum\limits_{\substack{i \in \mcl{N}_{w} : \discV_i = 1}}(1 - \discV'_i)$. Hence, if $\discV'_{i} = 1\text{ for all } i \in  \left\{ \mcl{N}_w : \discV_i = 1\right\}$, then $\sum\limits_{\substack{i \in \mcl{N}_{w} : \discV_i = 1}}(1 - \discV'_i) = 0$ and $\phi(\discV, \discV') = \rho(\discV, \discV')$. In all the other cases, $\sum\limits_{\substack{i \in \mcl{N}_{w} : \discV_i = 1}}(1 - \discV'_i) > 0$ and therefore $\rho(\discV, \discV') < \phi(\discV, \discV')$, as desired. \Halmos\endproof\vspace{5pt}

\noindent\textsc{Proof of Theorem~\ref{th_ccg_convergence}.}
Since $\ySet$ is finite, in the worst-case $|\ySet|$ steps are needed. In fact, after $|\ySet|$ steps, set $\widehat{\ySet}$ will be equivalent to set $\ySet$ since in each iteration (\ref{subproblem_ccg}) generates a new variable $\yvk{} \notin \widehat{\ySet}$ if $ub - lb > 0$.\Halmos\endproof

\section{Extensions of the Exact Algorithm}\label{sec:exa_extension}

Algorithms~\ref{Algo_Logic_Benders}~and~\ref{Algo_CCG} can be further extended or slightly modified to tackle a richer class of problems. Below we list some of the possible extensions.

\noindent\textbf{Finite $\Xi$ and mixed integer $\mathcal{Y}$:} If problem~\eqref{exact_phi_definition} involves a finite uncertainty set $\Xi$ and mixed integer set $\mathcal{Y}$,  Algorithm~\ref{Algo_CCG} will produce the optimal solution in no more than $|\Xi|$ iterations. Intuitively, this is the case as there are $|\Xi|$ inner $\displaystyle\min_{\bm y\in\mathcal{Y}} \max_{\bm \xi\in\Xi(\bm w,\bm{\bar \xi} )} \; \bm \xi^\top \bm C \bm w  + \bm \xi^\top \bm P \bm y$ problems, one for each  $\bm{\bar \xi}\in\Xi$. Therefore, in the worst-case the column-and-constraint generation algorithm will have to evaluate the worst-case cost of all $|\Xi|$ values of $\bm{\bar \xi}$ to find the optimal solution.

\noindent\textbf{Continuous $\Xi$  and mixed integer $\mathcal Y$:} If problem~\eqref{exact_phi_definition} involves a continuous uncertainty set $\Xi$ and mixed integer set $\mathcal Y$,  Algorithm~\ref{Algo_CCG} will produce the optimal solution in a finite number of iterations. In particular, consider the case where $\mathcal{Y} := \{\bm y_B\in\{0,1\}^{N_{\bm y_B}}, \bm y_C\in\mathbb{R}^{N_{\bm y_C}}: \bm F \bm y_B +  \bm G \bm y_C \leq \bm h\}$ is a bounded set for some matrices $\bm F\in \mathbb{R}^{m\times N_{\bm y_B}},\,\bm G\in \mathbb{R}^{m\times N_{y_C}}$ and vector $\bm h\in \mathbb{R}^{m}$. We can define sets $\mathcal{Y}_B:= \{\bm y_B\in\{0,1\}^{N_{\bm y_B}}: \exists \bm y_C\in\mathbb{R}^{N_{\bm y_C}} \text{ s.t. }  \bm F \bm y_B +  \bm G \bm y_C \leq \bm h\}$ and $\mathcal{Y}_C(\bm y_B):= \{\bm y_C\in\mathbb{R}^{N_{\bm y_C}}: \bm F \bm y_B +  \bm G \bm y_C \leq \bm h\}$. Problem~\eqref{exact_phi_definition} can be written as
\begin{equation*}
\begin{array}{ll}
\Phi(\discV) & = \displaystyle\max_{\bm{\bar \xi}\in\Xi} \min_{\bm y_B\in \mathcal{Y}_B}\min_{\bm y_C\in \mathcal{Y}_C(\bm y_B)} \max_{\bm \xi\in\Xi(\bm w, \bm{\bar \xi})} \; \bm \xi^\top \bm C \bm w  + \bm \xi^\top \bm P_B \bm y_B + \bm \xi^\top \bm P_C \bm y_c\\

& = \displaystyle\max_{\bm{\bar \xi}\in\Xi} \min_{\bm y_B\in \mathcal{Y}_B} \max_{\bm \xi\in\Xi(\bm w, \bm{\bar \xi})}\min_{\bm y_C\in \mathcal{Y}_C(\bm y_I)} \; \bm \xi^\top \bm C \bm w  + \bm \xi^\top \bm P_B \bm y_B + \bm \xi^\top \bm P_C \bm y_c\\

& = \displaystyle\max_{\bm{\bar \xi}\in\Xi} \min_{\bm y_B\in \mathcal{Y}_B} \max_{\bm \xi\in\Xi(\bm w, \bm{\bar \xi})}\max_{\bm \lambda\in\Lambda(\bm \xi)} \; \bm \xi^\top \bm C \bm w  + \bm \xi^\top \bm P_B \bm y_B + \bm \lambda^\top \bm F \bm y_B - \bm \lambda^\top\bm h,\\
\end{array}
\end{equation*}
where $\Lambda(\bm \xi) = \{\bm \lambda\in\mathbb{R}^m_+: (\bm P_C)^\top \bm \xi + \bm G^\top \bm \lambda = 0\}$. The second equality follows from the min-max theorem since both sets $\Xi(\bm w,\bm{\bar\xi})$ and $\mathcal{Y}_C(\bm y_B)$ are convex and $\mathcal{Y}_C(\bm y_B)$ is bounded by construction, while the last equality results from the dualization of $\displaystyle\min_{\bm y_C\in \mathcal{Y}_C(\bm y_B)}\bm \xi^\top \bm P_C \bm y_C$. The last problem has a similar structure as \eqref{exact_phi_definition} and Algorithm~\ref{Algo_CCG} can be directly applied, producing the optimal solution in no more than $|\mathcal{Y}_B|$ iterations.

\noindent\textbf{Binary first stage decisions $\bm x\in\mathcal{X}$:} Binary first stage decisions $\bm x\in\mathcal{X}\subseteq \{0,1\}^{N_x}$ can be incorporated into the problem in two ways which will lead to two variants of the algorithm. The first way is to express the problem as
\begin{subequations}\label{binaryx1}
\begin{equation}\label{outerw}
\min_{\bm w \in \mathcal{W},\,\bm x\in\mathcal{X}}\;\Phi(\bm w,\bm x),
\end{equation}
where
\begin{equation}\label{inner}
    \Phi(\bm w,\bm x) =  \max_{\bm{\bar\xi}\in\Xi} \;\min_{\bm y\in\mathcal{Y}(\bm  x)}\; \max_{\bm\xi\in\Xi(\bm w,\bm{\bar\xi})} \;\bm\xi^\top \bm Q\bm x + \bm\xi^\top \bm C\bm w + \bm\xi^\top \bm P\bm y. 
\end{equation}
\end{subequations}
Since both $\bm w$ and $\bm x$ are binary vectors, Algorithm~\ref{Algo_Logic_Benders} can still be applied in conjunction with the Logic Benders inequalities~\eqref{exact_convex_envelope}. In this case, in the worst-case Algorithm~\ref{Algo_Logic_Benders} will converge in $|\mathcal{W}||\mathcal{X}|$ iterations. Notice that the feasible region of the recourse decisions $\mathcal{Y}(\bm x)$ can be affected by the first stage decisions. However, this does not impact Algorithm~\ref{Algo_CCG} which for each $\bm x$ will produce the optimal solution in $|\mathcal{Y}(\bm x)|$ iterations in the worst-case. If for a given $\bm x$ the corresponding feasible region $\mathcal{Y}(\bm x)$ results being empty, a feasibility cut must be added to problem~\eqref{outerw} in order to exclude the solution from the problem. Note that due to the binary nature of variables $\discV$, this can be easily extended to the case in which $\mathcal{Y}(\discV,\bm{x})$.

An alternative modelling approach is to include the first stage decision in the definition of $\Phi(\discV)$ and express the problem as 
\begin{subequations}\label{binaryx2}
\begin{equation}\label{outerw2}
\min_{\bm w \in \mathcal{W}}\;\Phi(\bm w),
\end{equation}
where
\begin{equation}\label{inner2}
    \Phi(\bm w) = \min_{\bm x\in\mathcal{X}}\; \max_{\bm{\bar\xi}\in\Xi} \;\min_{\bm y\in\mathcal{Y}(\bm  x)}\; \max_{\bm\xi\in\Xi(\bm w,\bm{\bar\xi})} \;\bm\xi^\top \bm Q\bm x + \bm\xi^\top \bm C\bm w + \bm\xi^\top \bm P\bm y. 
\end{equation}
\end{subequations}
In this setting, evaluating $\Phi(\bm w)$ requires to have an outer column-and-constraint generation algorithm to optimize over $\bm x$. This additional algorithm will iteratively  solve the following relaxation of problem \eqref{inner2} 
\begin{equation*}\label{MasterOfInner}
\begin{array}{rll}
    \min&\tau \\
    \text{s.t.} & \bm x\in\mathcal{X},\,\tau \in\mathbb{R},\,\bm y:\widehat\Xi\mapsto \mathbb{R}^{N_{\yvk{}}}\\
    &\displaystyle \max_{\bm\xi\in\Xi(\bm w,\bm{\bar \xi})} \bm\xi^\top \bm Q\bm x + \bm\xi^\top \bm C\bm w + \bm\xi^\top \bm P\bm y(\bm{\bar \xi}) \leq \tau&  \quad\forall\bm{\bar\xi}\in\widehat\Xi\\
    &\bm y(\bm{\bar \xi}) \in\mathcal{Y}(\bm x)& \quad\forall \bm{\bar\xi}\in\widehat\Xi,\\
    \end{array}
\end{equation*}
where $\widehat \Xi$ is a finite subset of $\Xi$, to obtain a progressive solution $\bm x^*$ and a lower bound on $\Phi(\bm w)$. Using $\bm x^*$ we subsequently evaluate $\displaystyle \max_{\bm{\bar\xi}\in\Xi} \;\min_{\bm y\in\mathcal{Y}(\bm  x^*)}\; \max_{\bm\xi\in\Xi(\bm w,\bm{\bar\xi})} \;\bm\xi^\top \bm Q\bm x^* + \bm\xi^\top \bm C\bm w + \bm\xi^\top \bm P\bm y $ using Algorithm~\ref{Algo_CCG} to generate an upper bound on the value of $\Phi(\bm w)$ and the next $\bm{\bar\xi}$ to be added in set $\widehat \Xi$. The algorithm terminates when the values of the  upper and lower bounds coincide. Similar to the previous algorithms, in the worst-case, this extra algorithmic layer will produce the optimal solution in $|\mathcal{X}|$ iterations.

\noindent\textbf{Mixed integer first stage decisions $\bm x\in\mathcal{X}$:} In the case where the first stage decisions $\bm x = (\bm x_B, \bm x_C)$ involve both binary $\bm x_B$ and continuous $\bm x_C$ decisions,  $\Phi(\bm w,\bm x_B, \cdot)$ in problem~\eqref{binaryx1} will in general be non-convex. Since Algorithm~\ref{Algo_Logic_Benders} approximates $\Phi(\cdot)$ using a convex piecewise linear underestimator, it  cannot be used to approximate  $\Phi(\bm w,\bm x_B, \cdot)$. We can, however, use formulation \eqref{binaryx2} and the column-and-constraint generation algorithm outlined above. Unfortunately, in the worst-case the algorithm will convergence to the optimal solution asymptotically, unlike when $\bm x$ is purely binary where we have finite worst-case convergence.

\noindent\textbf{Stochastic and distributionally robust optimization problems:} It is interesting to note that Algorithm~\ref{Algo_Logic_Benders} can be used for any two-stage optimization problem (stochastic or robust) with decision-dependent information discovery. In all cases, $\discV$ controls the timing of information and $\Phi(\bm w)$ represents the value of the two-stage problem for the information structure prescribed by $\bm w$. The only requirements are  $(i)$ to be able to exactly evaluate $\Phi(\bm w)$  so as Algorithm~\ref{Algo_Logic_Benders} can produce the optimal solution, and $(ii)$ perform the evaluation efficiently for the algorithm to be practical.

\section{Proofs for Section~\ref{sec:k_adapt}}

\noindent \textsc{Proof of Theorem~\ref{th_symmetry_breaking}.}
We show that for each feasible solution of problem (\ref{mod_bi_problem}) an equivalent solution in which $\alDV$ belongs to set  $\mcl{A}_S$ exists.
Let $\sol$ be a feasible solution of problem (\ref{mod_bi_problem}) such that for 
some $k \in \Kset \setminus\{K\}$, $\alDV_{k} < \alDV_{k+1}$. We prove that an 
equivalent solution such that $\alDV_{k} \geq \alDV_{k+1}$ always exists, and show how to generate such a solution.
If the values of variables $\left(\alDV_k,\ \betakDV{k},\ \gammaDV{k},\ \yvk{k}\right)$ are interchanged with the values of variables $\left(\alDV_{k+1},\ \betakDV{k+1}\, \gammaDV{k+1}, \yvk{k+1}\right)$, we obtain a solution that satisfy (\ref{symmetry_breaking}). We need now to show that this solution is feasible if $\sol$ is feasible and that has the same cost. Obviously, the solution satisfy constraints (\ref{mod_bil:cnstr:balance_k}) and (\ref{mod_bil:cnstr:balance}). Furthermore, the values of variables $\betaDV$ remain unchanged, and hence the objective value does not change as well. We conclude that the two solutions are equivalent, as desired. \Halmos\endproof

\noindent \textsc{Proof of Theorem~\ref{th_alpha_bounds}.}
\label{proof_alpha_bound}We prove the theorem by demonstrating that 
$\mathcal{A}_S\subseteq\mathcal{A}_B$, thus leveraging on Theorem~\ref{th_symmetry_breaking} to 
demonstrate that at least one optimal solution of problem~\eqref{mod_bi_problem} satisfies 
$\alDV\in\mathcal{A}_B$.
For $K = 1$ the statement is trivial. Consider some $K>1$. We prove the statement by contradiction. Consider $k = 1$ and assume that there exists $\alDV\in\mathcal{A}_S$ with $\alDV_1<\frac{1}{K}$. Since $\alDV_1 > 0$, there exists $\epsilon \in (0,\frac{1}{K})$
such that $\alDV_1 =\frac{1}{K} - \epsilon$. Constraint $\bm e^\top \alDV = 1$ implies that
\begin{subequations}
\begin{equation}\label{alpha_bound_eq1}
    \sum_{k \in \Kset \setminus \{1\}}\alDV_{k} = 1 - \frac{1}{K} + \epsilon,
\end{equation}
while constraints $\alDV_k \geq \alDV_{k+1}$ for all $k \in \Kset \setminus \left\{K\right\}$ imply
\begin{equation}\label{alpha_bound_eq2}
   \alDV_{1} \geq \max_{k\in \Kset \setminus \{1\}}\alDV_{k}.
\end{equation}
We show that none of the $\alDV \in \rrr^{K}_{+}$ such that \eqref{alpha_bound_eq1} is satisfied satisfies \eqref{alpha_bound_eq2}, since
\begin{equation}
\alDV_1  \geq \min_{\alDV_{-1} \geq 0: \eqref{alpha_bound_eq1}}\; \max_{k \in \Kset \setminus\{1\}} \alDV_k =  \frac{1 + \epsilon - \frac{1}{K}}{K-1} \implies
\epsilon \leq 0.
\end{equation}
\end{subequations}
Therefore, we reach the contradiction that $\epsilon \in (0,\frac{1}{K})$. 
Consider now the generic component $k > 1$ and assume that $\alDV_k >\frac{1}{k}$. This implies that there exists $\epsilon \in \left(0, 1- \frac{1}{k}\right)$ such that $\alDV_k = \frac{1}{k} + \epsilon$. Constraint $\bm{e}^\top \alDV = 1$ implies
\begin{subequations}
\begin{equation}\label{alpha_bound_eq3}
\sum_{\til{k} \in \Kset \setminus \{k\}}\alDV_{\til{k}} = 1 - \frac{1}{k} - \epsilon,
\end{equation}
while constraint  constraints $\alDV_k \geq \alDV_{k+1}$ for all $k \in \Kset \setminus\left\{K\right\}$ imply
\begin{equation}\label{alpha_bound_eq4}
\min_{ \left\{ \til{k} \in \Kset : \til{k} < k\right\} }\alDV_{\til{k}} \geq \alDV_{k} \geq \max_{ \left\{ \til{k} \in \Kset : \til{k} > k\right\} }\alDV_{\til{k}}.
\end{equation}
\end{subequations}
We show that none of the $\alDV$ that satisfies \eqref{alpha_bound_eq3} will satisfy \eqref{alpha_bound_eq4}, since
\begin{equation*}
\max_{\alDV_{-k} \geq 0 :\eqref{alpha_bound_eq3}}\; \min_{\tilde 
k \in\{ \Kset : \tilde k < k\}}\alDV_{\til{k}} \geq \alDV_k\implies
\frac{1}{k} + \epsilon  \leq \frac{1 - \epsilon - \frac{1}{k}}{k-1} \implies
\epsilon \leq 0.
\label{alpha_bound_kl_eq_end}
\end{equation*}
Therefore, we reach the contradiction that $\epsilon \in (0,1-\frac{1}{k})$. \Halmos\endproof

\noindent \textsc{Proof of Theorem~\ref{th_optimistic_cut}.}
Let $(\bm w, \bm \alpha, \bm \beta,\{\bm y^k,\bm \beta^k,\bm \gamma^k\}_{k\in\mathcal{K}})$ be an optimal solution to problem~\eqref{mod_bi_problem_lin}. For all $k'\in\mathcal{K}$ for which $\bm\alpha_{k'} = 0$, we construct the solution where we set $\bm{y}^{k'} := \bm y^k$ for some $k\in\mathcal{K}$ whose corresponding $\bm{\alpha}_k>0$, which from Observation~\ref{observation} we know is also an optimal solution. Using this modified optimal solution, we have that
\begin{subequations}
\begin{align}
         &  \displaystyle\max_{\bm{\bar\xi}\in\Xi} \;\min_{k\in\mathcal{K}}\;\max_{\bm\xi\in\Xi(\bm w,\bm{\bar\xi})} \; \bm \xi^\top\bm C\bm w + \bm\xi^\top \bm P\bm y^k\nonumber\\
         =&\displaystyle\max_{\bm{\bar\xi}\in\Xi} \;\max_{\bm\xi^k\in\Xi(\bm w,\bm{\bar\xi})}\;\min_{k\in\mathcal{K}} \; {\bm \xi^k}^\top\bm C\bm w + {\bm \xi^k}^\top \bm P\bm y^k\label{optimize_over_xi}\\
         =&\displaystyle\min_{k\in\mathcal{K}} \; {\bm \xi^{k*}}^\top\bm C\bm w + {\bm \xi^{k*}}^\top \bm P\bm y^k\label{min_k}\\
      =&\,\displaystyle {\bm \xi^{k*}}^\top\bm C\bm w + {\bm \xi^{k*}}^\top \bm P\bm y^k\quad \forall k\in\mathcal{K}\label{equality}\\
        \geq&\,\displaystyle {\bm \zeta^{\bm w}}^\top \bm w + {\bm \zeta^{\bm y}}^\top \bm y^k \quad \forall k\in\mathcal{K},
\end{align}
\end{subequations}
where \eqref{optimize_over_xi}  follows from \cite[Lemma 1]{vayanos2020robust}. Vectors $\{\bm \xi^{k*}\}_{k\in\mathcal{K}}$ indicate the optimal solution of \eqref{optimize_over_xi}, which we substitute resulting in \eqref{min_k}. Equality \eqref{equality} follows since by construction of the optimal solution ${\bm\xi^{k*}}^\top\bm C\bm w + {\bm \xi^{k*}}^\top \bm P\bm y^k$ is equal for all $k\in\mathcal{K}$, while the last inequality results from the definition of $\bm \zeta^{\bm w}$ and $\bm \zeta^{\bm y}$. From strong duality we also have that
\begin{equation*}
    \begin{array}{ll}
         &  \displaystyle\max_{\bm{\bar\xi}\in\Xi} \;\min_{k\in\mathcal{K}}\;\max_{\bm\xi\in\Xi(\bm w,\bm{\bar\xi})} \; \bm \xi^\top\bm C\bm w + \bm\xi^\top \bm P\bm y^k\\[2ex]
         =&\left[
         \begin{array}{ll}
    \begin{array}{cl}
         \min &\displaystyle \quad   {\bm b}^\top\left( {\bm \beta}  + \sum_{k\in \mathcal K}  {\bm \beta}^k\right)\\
         \text{s.t.} & 
          \quad {\bm \alpha} \in \mathbb R^K_+, \; {\bm \beta} \in \mathbb R^L_+, \; {\bm \beta}^k \in \mathbb R^L_+, \; {\bm \gamma}^k \in \mathbb R^{N_\xi}, \; k \in \mathcal K \\
        & \quad \bm e^\top {\bm \alpha}  = 1  \\
        & \quad  {\bm A}^\top{\bm \beta}  =  \displaystyle \sum_{k \in \mathcal K} {\bm w} \circ  {\bm \gamma}^k \\
         & \quad     {\bm A}^\top{\bm \beta}^k +  {\bm w} \circ {\bm \gamma}^k = {\bm \alpha}_k \left( {\bm C} {\bm w}  + {\bm P} {\bm y}^k \right) \quad \forall k \in \mathcal K \\
    \end{array}
         \end{array}
         \right] 
    \end{array}
\end{equation*}
hence the result follows.
\Halmos\endproof

\newpage
\section{Supporting Tables for Numerical Examples and Case Study}

\begin{table}[htbp!]
\begin{footnotesize}
  \begin{center}
  \caption{Instance description}
          \sffamily
\renewcommand{\arraystretch}{0.8}
\begin{tabular}{cc}
    \begin{tabular}{|c|c|l|}
    \hline
    \hline
    Network & $U$  &  \multicolumn{1}{c|}{$T$}\\
    \hline
    TS1N15  &0.10 &\makecell[l]{5, 10, 15, 20, 25, 30, 35,\\ 40, 45, 50, 55, 60, 65, 70} \\
    \hline
    TS2N10  &0.20 & 15, 20, 23, 25, 27, 30, 32, 35, 38\\
    \hline
    TS3N16  &0.10 &\makecell[l]{15, 20, 25, 30, 35, 40, 45, 50, 55,\\ 60, 65, 70, 75, 80}\\
     \Xhline{1pt}
    \end{tabular}
      & \hfill 
       \begin{tabular}{|c|c|l|}
    \hline
    \hline
    Network & $U$  &  \multicolumn{1}{c|}{$T$}\\
    \hline
     TS1N30  &0.05 & \makecell[l]{5, 10, 15, 20, 25, 30, 35, 40, 46, 50,\\ 55, 60, 65, 70, 73, 75, 80, 85} \\
     \hline
     TS2N19  &0.15 & 15, 20, 23, 25, 27, 30, 32, 35, 38, 40, 45\\
     \hline
     TS3N31 & 0.05 &\makecell[l]{ 15, 20, 25, 30, 35, 40, 45, 50, 55, 60, 65,\\ 70, 75, 80, 85, 90, 110, 105, 110}\\
     \Xhline{1pt}
    \end{tabular}
    \end{tabular}\label{tab:op_instances}%
    \end{center}
    \end{footnotesize}
    \end{table}

\begin{table}[!htbp]
  \centering
  \caption{Location areas and departments of the Alrijne hospital clustered based on the proximity of the departments on each floor. Column Node and Department denote the cluster in which each department is assigned, for example node 1A covers Nuclear Medicine and Radiology. The column Expected Crates reports the expected percentages of crates at each node derived from data used in \cite{alrijneBachelorThesis}.}
      \begin{scriptsize}
        \sffamily
\setlength{\tabcolsep}{5pt}
\renewcommand{\arraystretch}{1.3}
    \begin{tabular}{|clc|clc|}
    \hline
    \hline
   Node & \multicolumn{1}{l}{Department} &\makecell{Expected \\Crates} &\multicolumn{1}{l}{Node} & Department &\makecell{Expected \\Crates}\\
    \hline
    
    0B & \multicolumn{1}{l}{Gynecology outpatient clinic} & 5\% & 2A    &  Urology ward &10\%  \\
    
   \multirow[t]{5}[0]{*}{0C} & Orthopedist &10\%  & \multirow[t]{2}[0]{*}{2B} &  Day care/short stay \& orthopedics & 5\% \\
    
    & Ophthalmology outpatient clinic & &      &  Surgery ward & \\
          
     & \multicolumn{1}{l}{Ear, nose and throat outpatient clinic} & & \multirow[t]{2}[0]{*}{2C} &  Anesthesia & 6\%\\
          
    & Orthopedics outpatient clinic &    &   &  Operating theatres &\\
          
    & Wound center && \multirow[t]{2}[0]{*}{3A} &  Geriatric trauma unit & 7\% \\
          
  \multirow[t]{2}[0]{*}{1A} & Nuclear medicine & 10\% &     &   Neurology  &  \\
          
    & Radiology & & \multirow[t]{2}[0]{*}{3B} &  Cardiac care unit   &5\% \\
           
   1B & Poli plastic surgery &7\% &      &   Cardiology II &  \\
    
    \multirow[t]{4}[0]{*}{1C} & Cardiology I &  5\%& 4A    &Mother and child center & 7\%\\
    
    & \multicolumn{1}{l}{Lung diseases outpatient clinic}  &  & \multirow[t]{2}[0]{*}{4B}&  Digestive and liver diseases & 8\% \\
    
    & Endoscopy &   &    &   Oncology ward  & \\
    
    & Dermatology outpatient clinic &     &5A  &Child and youth department &10\%\\

    &   &    & 5B& Dialysis & 5\% \\

    \Xhline{1pt}
    \end{tabular}%
          \end{scriptsize}
  \label{tab:Alrijne:departments}%
\end{table}%

\end{document}